 \documentclass[10pt,twoside,english]{article}
\usepackage{amssymb,amsmath,babel,geometry,latexsym,graphics,tabularx,shapepar,enumerate}
\ifx\optionkeymacros\undefined\else \fi

\catcode`\Œ=\active\defŒ{{\aa}}       
\catcode`\º=\active\defº{\int}        
\catcode`\=\active\def{\c c}        
\catcode`\¶=\active\def¶{\partial}    
\catcode`\Ä=\active\defÄ{\oint}       
\catcode`\Æ=\active\defÆ{\triangle}   
\catcode`\Â=\active\defÂ{\neg}        
\catcode`\µ=\active\defµ{\mu}         
\catcode`\¿=\active\def¿{{\o}}        
\catcode`\¹=\active\def¹{\pi}         
\catcode`\Ï=\active\defÏ{{\oe}}       
\catcode`\§=\active\def§{{\ss}}       
\catcode`\ =\active\def {\dagger}     
\catcode`\Ã=\active\defÃ{\sqrt}       
\catcode`\·=\active\def·{\Sigma}      
\catcode`\Å=\active\defÅ{\approx}     
\catcode`\½=\active\def½{\Omega}      
\catcode`\£=\active\def£{{\it\$}}     
\catcode`\°=\active\def°{\infty}      
\catcode`\¤=\active\def¤{{\S}}        
\catcode`\¦=\active\def¦{{\P}}        
\catcode`\¥=\active\def¥{\bullet}     
\catcode`\»=\active\def»{\leavevmode\raise.585ex\hbox{\b a}}      
\catcode`\¼=\active\def¼{\leavevmode\raise.6ex\hbox{\b o}}        
\catcode`\­=\active\def­{\not=}       
\catcode`\²=\active\def²{\leq}        
\catcode`\³=\active\def³{\geq}        
\catcode`\Ö=\active\defÖ{\div}        
\catcode`\É=\active\defÉ{{\dots}}     
\catcode`\¾=\active\def¾{{\ae}}       
\catcode`\Ç=\active\defÇ{\ll}         
\catcode`\Ò=\active\defÒ{``}          
\catcode`\Á=\active\defÁ{!`}          
\catcode`\¢=\active\def¢{\rlap/c}     
\catcode`\Ô=\active\defÔ{`}           
\catcode`\Õ=\active\defÕ{'}           

\catcode`\=\active\def{{\AA}}       
\catcode`\'=\active\def'{\c C}        
\catcode`\¯=\active\def¯{{\O}}        
\catcode`\¸=\active\def¸{\Pi}         
\catcode`\Î=\active\defÎ{{\OE}}       
\catcode`\®=\active\def®{{\AE}}       
\catcode`\×=\active\def×{\diamond}    
\catcode`\¡=\active\def¡{\accent'27}  
\catcode`\Ó=\active\defÓ{''}          
\catcode`\±=\active\def±{\pm}         
\catcode`\È=\active\defÈ{\gg}         
\catcode`\À=\active\defÀ{?`}          
\catcode`\Ð=\active\defÐ{--}          
\catcode`\Ñ=\active\defÑ{---}         

\catcode`\Š=\active\defŠ{\"a}        
\catcode`\'=\active\def'{\"e}        
\catcode`\•=\active\def•{\"{\i}}     
\catcode`\š=\active\defš{\"o}        
\catcode`\Ÿ=\active\defŸ{\"u}        
\catcode`\Ø=\active\defØ{\"y}        
\catcode`\€=\active\def€{\"A}        
\catcode`\…=\active\def…{\"O}        
\catcode`\†=\active\def†{\"U}        
\catcode`\‡=\active\def‡{\'a}        
\catcode`\Ž=\active\defŽ{\'e}        
\catcode`\'=\active\def'{\'{\i}}     
\catcode`\—=\active\def—{\'o}        
\catcode`\œ=\active\defœ{\'u}        
\catcode`\ƒ=\active\defƒ{\'E}        
\catcode`\ˆ=\active\defˆ{\`a}        
\catcode`\=\active\def{\`e}        
\catcode`\"=\active\def"{\`{\i}}     
\catcode`\˜=\active\def˜{\`o}        
\catcode`\=\active\def{\`u}        
\catcode`\Ë=\active\defË{\`A}        
\catcode`\‹=\active\def‹{\~a}        
\catcode`\–=\active\def–{\~n}        
\catcode`\›=\active\def›{\~o}        
\catcode`\Ì=\active\defÌ{\~A}        
\catcode`\"=\active\def"{\~N}        
\catcode`\Í=\active\defÍ{\~O}        
\catcode`\‰=\active\def‰{\^a}        
\catcode`\=\active\def{\^e}        
\catcode`\"=\active\def"{\^{\i}}     
\catcode`\™=\active\def™{\^o}        
\catcode`\ž=\active\defž{\^u}        

\let\optionkeymacros\null

\usepackage{hyperref}
\usepackage[all,2cell]{xy} \UseAllTwocells \SilentMatrices

\newtheorem{Theorem}{Theorem}
\newtheorem{Lemma}[Theorem]{Lemma}
\newtheorem{Proposition}[Theorem]{Proposition}
\newtheorem{Conjecture}[Theorem]{Conjecture}
\newtheorem{Corollary}[Theorem]{Corollary}
\newtheorem{Definition}[Theorem]{Definition}
\newtheorem{Remark}[Theorem]{Remark}

\def\pp{\mathfrak{p}}
\def\aa{\mathfrak{a}}
\def\SS{H}
\def\cc{\chi}

\newcommand{\ZZ}{\mathbb{Z}}
\newcommand{\FF}{\mathbb{F}}
\newcommand{\CC}{\mathbb{C}}
\newcommand{\QQ}{\mathbb{Q}}

\newcommand{\LL}{\mathbb{L}}

\newcommand{\TT}{\mathbb{T}}

\newcommand\CVD{{\hfill\hfil{\lower 2 pt\hbox{\vrule\vbox to 7pt 
{\hrule width 6pt\vfill\hrule}\vrule}}}\vskip 0.5cm}

\title{Universal Gauss-Thakur sums and $L$-series
\;\footnote{The second author was supported by the contract ANR ``HAMOT", BLAN-0115-01}}
\author{
Bruno Angl\`es\footnote{Current address: LMNO, UniversitŽ de Caen BP 5186
F 14032 Caen Cedex.} \; \& 
Federico Pellarin\footnote{Current address: ICJ Lyon and Saint-Etienne, 23, rue du Dr. Paul Michelon, 42023 Saint-Etienne Cedex.} }

\begin{document}

\maketitle

\begin{small}
\noindent\textbf{Abstract.} In this paper we study the behavior of the function $\omega$ of Anderson-Thakur  (introduced in \cite{AT})
evaluated at the elements of the algebraic closure $\FF_q^{\text{alg}}$ of the finite field with $q$ elements $\FF_q$. Indeed, this function has quite a remarkable connection with explicit class field theory 
for the field $K=\FF_q(\theta)$. We will see that these values, together with the values at $\FF_q^{\text{alg}}$
of its divided derivatives, generate over $\FF_q^{\text{alg}}$ the maximal abelian extension of $K$ which is tamely 
ramified at infinity. We will also see that $\omega$ is, in a way that we will explain in detail, an {\em universal Gauss-Thakur} sum. We will then use 
these results to show the existence of functional relations for a class of $L$-series introduced by the second author in \cite{Pe}.
Our results will be finally applied to obtain a new class of congruences for Bernoulli-Carlitz fractions, and an analytic conjecture 
is stated, implying an interesting behavior of such fractions modulo prime ideals of $A=\FF_q[\theta]$.
\end{small}

\medskip

\tableofcontents

\section{Introduction, results}

The present paper is divided in two parts (Section \ref{kroneckerweber} 
for the first part and Sections \ref{sectionthree} and \ref{congruences} for the second), both motivated by the interesting behavior that the function $\omega$ of {\em Anderson and Thakur} (\footnote{Introduced in \cite{AT}.}) exhibits at the roots of unity, and the consequence that this behavior has on analytic properties of certain $L$-series introduced in \cite{Pe}.

We will first be concerned with the values of the function $\omega$ at 
the roots of unity and we will prove, among several results, 
Theorem \ref{mtre}, which provides, we hope, an alternative approach to {\em explicit class field theory}.
We will also prove, in Theorem \ref{gaussthakur}, that $\omega$ is, in a certain sense that will be made more precise later, an {\em universal Gauss-Thakur sum}. 

Theorem \ref{gaussthakur} will be used in 
in the second part, where we will consider a class of $L$-series that was recently introduced by the second author in \cite{Pe} and we will study their behavior at the roots of unity. We will prove, in Theorem \ref{functionalrelations}, functional identities in the same vein as in loc. cit., but in a much more general 
multivariable setting. Among others, 
some applications to {\em Bernoulli-Carlitz} numbers are given, in Theorem \ref{bernoulli}. 

\medskip

Here is, more specifically, the content of the present paper.
Our purpose in Section \ref{kroneckerweber}, the first part of this paper, is to focus on explicit class field theory for the field $K$. The classical Kronecker-Weber theorem states that
the maximal abelian extension $\QQ^{\text{ab}}$ of the field of rational numbers $\QQ$ in the field of complex numbers $\CC$ is 
generated by the values of the exponential function $$e^z=\sum_{n\geq 0}\frac{z^n}{n!}$$ at the elements $\sqrt{-1}\pi\rho$, $\rho\in\QQ$, or, in other words,
by the complex roots of the polynomials
$$X^n-1,\quad n\geq 1.$$
The prominency of an analytic function in an algebraic problem is the essence of Kronecker's {\em Jugendtraum} (it later became the {\em twelfth Hilbert's problem}) 
and was confirmed in other situations by other authors, namely by Hayes in 1974, which in \cite{Ha1} analytically expressed a minimal set
of generators of the maximal
abelian extension of the field $K=\FF_q(\theta)$  tamely ramified at infinity by means of the torsion values of {\em Carlitz's exponential function}, and further constructed the maximal abelian 
extension of $K$ again relating it to the torsion of Carlitz's module.

We denote by $v_\infty$ the $\theta^{-1}$-adic valuation normalized, in all the following, by setting $v_\infty(\theta)=-1$.
Let $K_\infty$ be the completion of $K$ for $v_\infty$, and let us consider the completion $\CC_\infty$ of an algebraic closure of $K_\infty$
for the unique extension of this valuation, in which we embed an algebraic closure of $K$. Carlitz's exponential function is the surjective, $\FF_q$-linear, rigid analytic entire function
$$\exp:\CC_\infty\rightarrow\CC_\infty$$
defined by
$$\exp(z)=\sum_{n\geq0}\frac{z^{q^n}}{d_n},$$
where 
$$d_0=1,\quad d_n=(\theta^{q^n}-\theta)(\theta^{q^n}-\theta^q)\cdots(\theta^{q^n}-\theta^{q^{n-1}}),\quad n>0.$$
The kernel of this function turns out to be generated by a {\em period} $\widetilde{\pi}$, unique up to multiplication by an 
element of $\FF_q^\times$, that can be computed by using the following product expansion
\begin{equation}\label{pi}
\widetilde{\pi}:=\theta(-\theta)^{\frac{1}{q-1}}\prod_{i=1}^\infty(1-\theta^{1-q^i})^{-1}\in (-\theta)^{\frac{1}{q-1}}K_\infty,
\end{equation}
once a $(q-1)$-th root of $-\theta$ is chosen. 

\medskip

\noindent\emph{Anderson-Thakur function.} This function, introduced in \cite[Proof of Lemma 2.5.4 p. 177]{AT}, is defined by the infinite product 
\begin{equation}\label{definitionomega}
\omega(t)=(-\theta)^{\frac{1}{q-1}}\prod_{i\geq0}\left(1-\frac{t}{\theta^{q^{i}}}\right)^{-1}\end{equation} (it is customary to make the same choice as in (\ref{pi}) for the $(q-1)$-th root),
converges for $t\in\CC_\infty$ such that $|t|\leq 1$ (where $|\cdot|$ is an absolute value associated to 
$v_\infty$) and can be extended to a non-vanishing rigid analytic function over $$\CC_\infty\setminus\{\theta^{q^k};k\geq 0\}$$ with simple poles at $\theta^{q^k}$,
$k\geq 0$. In \cite{gamma}, many analogies with 
{\em Euler's gamma function} are tracked. For instance, variants of the {\em translation formula}, {\em Gauss multiplication formulas}
and {\em reflection formulas} for the gamma function hold for $\omega$. We are going to study yet another property of $\omega$.

For $\zeta\in \FF_q^{\text{alg}}$, the product $\omega(\zeta)$ in (\ref{definitionomega}) converges to an algebraic element 
of $\CC_\infty$.
More generally, consider the $\CC_\infty$-linear {\em divided derivatives} $$\mathcal{D}_n:\CC_\infty[[t]]\rightarrow\CC_\infty[[t]],\quad n\geq0$$
defined by setting $$\mathcal{D}_nt^m=\binom{m}{n}t^{n-m}.$$
Then, for all $n\geq 0$, the formal series $\mathcal{D}_n\omega$ converges for all $t\in\CC_\infty$ such that $|t|\leq 1$. For $\zeta\in\FF_q^{\text{alg}}$,
the series $(\mathcal{D}_n\omega)(\zeta)$ converges in fact to an element of $K^{\text{alg}}$ (these properties will be apparent in the paper).
Furthermore, let $E^\infty$ be the smallest subfield of $\CC_\infty$ containing the algebraic closure $\FF_q^{\text{alg}}$ of $\FF_q$ in $\CC_\infty$ and the values $$\exp(\widetilde{\pi}\rho),\quad \rho\in K.$$
The first result of this paper is the following.
\begin{Theorem}\label{mtre}
The field $E^\infty$ is also generated over $\FF_q^{\text{alg}}$
by the values $(\mathcal{D}_n\omega)(\zeta)$ for all $\zeta\in\FF_q^{\text{alg}}$ and $n\geq 0$.
\end{Theorem}
According to Hayes \cite{Ha1}, $E^\infty$ is equal to the maximal abelian extension of $K$ tamely ramified at $\infty$ 
in $\CC_\infty$ (\footnote{That is, the maximal abelian extension of $K$ whose perfection is contained in the subfield of Newton-Puiseux series $\cup_{n\geq 1}\FF_q^{\text{alg}}((\theta^{-1/n}))$.}), which obviously yields the next Corollary.
\begin{Corollary}
The higher derivatives of the function $\omega$ evaluated
at the elements of $\FF_q^{\text{alg}}$ generate, over $\FF_q^{\text{alg}}$, the maximal extension 
of $K$ which is abelian and tamely ramified at the infinity place.
\end{Corollary}

In the proof of Theorem \ref{mtre}, new functions generalizing the function $\omega$ are introduced. These are the functions
$\omega_{\mathfrak{a},j}$ of Subsection \ref{omegaaj}. 
They generalize the function $\omega=\omega_{\theta,0}$ in the sense of Proposition \ref{basis} and should be considered of same relevance
as $\omega$ itself, being associated to the kernel of $\phi_\mathfrak{a}$ (the image of a monic polynomial $\mathfrak{a}\in A$ by Carlitz's module)
in the same way as $\omega$ is associated to the kernel of $\phi_\theta$. These features 
will be discussed in detail in Subsection \ref{omegaaj}.

The second result of Section \ref{kroneckerweber}, closely related to Theorem \ref{mtre}, draws a portrait of Anderson-Thakur's function itself, as an {\em universal Gauss-Thakur sum}. This
analogue of Gauss sums, in $K^{\text{ab}}$, was introduced by Thakur in \cite{Tha1}. Thakur established
several analogues of classical results about Gauss sums such as Stickelberger factorization theorem and Gross-Koblitz formulas and other 
analogues of classical results
(see for example \cite{Tha1,Tha2,Tha3}). We refer to 
Subsection \ref{gaussthakursubs} for the background on Gauss-Thakur sums. We are going to describe a direct connection between Gauss-Thakur sums and the function $\omega$. 

Let $\pp$ be an irreducible monic polynomial of $A$ of degree $d$, let
$\Delta_\pp$ be the Galois group of the $\pp$-cyclotomic function field extension $K(\lambda_\pp)$ of $K$, where $\lambda_\pp$
is a non zero $\pp$-torsion element of $K^{\text{alg}}$. 
Gauss-Thakur sums can be associated to the elements of the dual character group $\widehat{\Delta}_\pp$ via the Artin symbol (see \cite[Sections 7.5.5 and 9.8]{Go2}). If $\chi$ is in $\widehat{\Delta}_\pp$, we denote by $g(\chi)$ the associated
Gauss-Thakur sum. In particular, we have the element $\vartheta_\pp\in\widehat{\Delta}_\pp$ obtained by reduction of the {\em Teichm\"uller character} \cite[Definition 8.11.2]{Go2}, uniquely determined by a choice of a root $\zeta$ of $\pp$, and the Gauss-Thakur sums $g(\vartheta_\pp^{q^j})$ associated to its $q^j$-th powers, with $j=0,\ldots,d-1$, which can be considered as
the building blocks of the Gauss-Thakur sums $g(\chi)$ for general $\chi\in\widehat{\Delta}_\pp$.

\begin{Theorem}\label{gaussthakur}
Let $\pp $ be a prime element of $A$ of degree $d$ and $\zeta$ a root of $\pp$ as above. We have:
$$g(\vartheta_\pp^{q^j})={\pp '(\zeta)}^{-q^j}\omega(\zeta^{q^j}),\quad j=0,\ldots,d-1.$$
\end{Theorem}
In this theorem, $\pp '$ denotes the derivative of $\pp $ with respect to $\theta$. We anticipate that 
Theorem \ref{gaussthakur} will play an important role in the proof of the next Theorem \ref{functionalrelations}. Also,
The Theorems \ref{mtre} and \ref{gaussthakur} are closely related. We will see, by Corollary \ref{corollary13} later in this paper,
that the field generated over $\FF_q^\text{alg}(\theta)$ by the various Gauss-Thakur sums $g(\vartheta_\pp)$, is also
equal to the field generated over $\FF_q^{\text{alg}}$ by the elements $\lambda\in K$ which are
$a$-torsion for $a\in A$ squarefree. But the proof of Theorem \ref{mtre} that we furnish, founded on an 
analytic formula (Proposition \ref{anidentity}), also tells us that
the last field is generated, over $\FF_q^{\text{alg}}$, by the elements $\omega(\zeta)$, $\zeta\in\FF_q^{\text{alg}}$
(see Proposition \ref{intermediateresult}).

\medskip

In Section \ref{sectionthree} we keep studying the values of $\omega$ at the elements of $\FF_q^\text{alg}$, but
we change our point of view by focusing now on certain $L$-series introduced in \cite{Pe}. Let $t$ be a variable in $\CC_\infty$
and let us consider the ring homomorphism $$\cc_{t}:A\rightarrow\FF_q[t]$$
defined by the formal replacement of $\theta$ by $t$. In other words, 
$\chi_t$ may be viewed as the unique ring homomorphism from $A$ to the ring of rigid analytic functions $\CC_\infty\rightarrow\CC_\infty$
such that $\chi_t(\theta)=t$. More generally, we shall consider $s$ independent variables $t_1,\ldots,t_s$ and consider 
the ring homomorphisms $$\chi_{t_i}:A\rightarrow\FF_q[t_1,\ldots,t_s],\quad i=1,\ldots,s$$ 
defined respectively by $\chi_{t_i}(\theta)=t_i$. To simplify our notations, we will write $\chi_{\xi}(a)$ or $a(\xi)$ for the evaluation at $t=\xi$ of the 
polynomial function $\chi_t(a)$ at a given element $\xi\in\CC_\infty$.
Let $\alpha$ be a positive integer and let $\beta_1,\ldots \beta_s$ be non-negative integers. The following formal series was introduced in \cite{Pe}: 
\begin{equation}\label{severalvar}L(\chi_{t_1}^{\beta_1}\cdots\chi_{t_s}^{\beta_s},\alpha)=\sum_{d\geq 0}\sum_{a\in A^+(d)}\chi_{t_1}(a)^{\beta_1}\cdots\chi_{t_s}(a)^{\beta_s}a^{-\alpha}\in
K_\infty[[t_1,\ldots,t_s]].\end{equation}
Here and in all the following, $A^+(d)$ denotes the set of monic polynomials of $A$ of degree $d$. It is easy to see that this series is well defined. 
As claimed in \cite[Remark 7]{Pe}, this series converges 
for all $(t_1,\ldots,t_s)\in\CC_\infty^s$ to a rigid analytic entire function of $s$ variables $t_1,\ldots,t_s$; see Proposition \ref{analytic}.

For the next result, we need further notation.
For $k$ a non-negative integer, we consider the $q$-ary expansion $k=k_0+k_1q+\cdots+k_sq^s$, where $k_0,k_1,\ldots,k_s$ are integers in the set
$\{0,\ldots,q-1\}$. We then denote by $\ell_q(k)$ the integer $k_0+k_1+\cdots+k_s$. The residue of $\omega(t)$ at $t=\theta$ 
is $-\widetilde{\pi}$:
$$\widetilde{\pi}=-\lim_{t\rightarrow\theta}(t-\theta)\omega(t).$$
In \cite[Theorem 1]{Pe}, it is proved that $$L(\chi_t,1)=\frac{\widetilde{\pi}}{(\theta-t)\omega(t)}.$$ Taking into account 
the functional equation 
$$\omega(t)^q=(t^q-\theta)\omega(t^q)$$ apparent in (\ref{definitionomega}), this implies that, for $m\geq 0$ integer,
$$V_{q^m,1}(t):=\widetilde{\pi}^{-q^m}L(\chi_t,q^m)\omega(t)=\frac{1}{(\theta^{q^m}-t)(\theta^{q^{m-1}}-t)\cdots(\theta-t)}.$$
This result provides an awaited connection between the function $\omega$ of Anderson and Thakur and the ``positive even" values of the {\em Goss zeta function} (or {\em Carlitz zeta values})
$$\zeta(n)=\frac{BC_n\widetilde{\pi}^n}{\Pi(n)},\quad n>0, \quad n\equiv0\pmod{q-1}$$
where $BC_n$ and $\Pi(n)$ denote respectively the $n$-th Bernoulli-Carlitz fraction and Carlitz's factorial of $n$, see Goss' book \cite[Section 9.1]{Go2}. Indeed, evaluating at $t=\theta$,
we get
$$L(\chi_\theta,q^m)=\zeta(q^m-1),\quad m\geq 1.$$

More generally, it is proved in \cite[Theorem 2]{Pe} that, if $\alpha\equiv1\pmod{q-1}$ and $\alpha\geq 1$,
then
$$\lambda_\alpha=\widetilde{\pi}^{-\alpha}L(\chi_t,\alpha)\omega(t)$$
is a rational function in $\FF_q(\theta,t)$. In \cite{Pe}, it is suggested that this result could be
a source of information in the study of the arithmetic properties of the Bernoulli-Carlitz fractions. 
However, the methods of loc. cit. (based on {\em deformations of vectorial modular forms} and Galois descent) are 
only partially explicit.

More recently, Perkins \cite{Per1} investigated the properties of certain {\em special polynomials} associated to variants of the functions 
$L(\chi_t^\beta,\alpha)$
with $\alpha\leq 0$ which turn out to be polynomial. He notably studied the growth of their degrees. Moreover, by using {\em Wagner's interpolation theory} for the map $\chi_t$, Perkins \cite{Per2} generalized some unpublished formulas of the second author and obtained explicit formulas for the series
$$L(\chi_{t_1}\cdots\chi_{t_s},\alpha),\quad \alpha>0,\quad 0\leq s\leq q,\quad \alpha\equiv s\pmod{q-1}.$$
We quote here a particular case of Perkins' formulas for the 
functions $L(\chi_t,\alpha)$ with $\alpha\equiv1\pmod{q-1}$:
\begin{equation}\label{perkinsformula}
L(\chi_t,\alpha)=\sum_{j=0}^{\mu}d_j^{-1}(t-\theta)(t-\theta^q)\cdots(t-\theta^{q^{j-1}})\zeta(\alpha-q^j)L(\chi_t,q^j),
\end{equation}
where $\mu$ is the biggest integer such that $q^\mu\leq \alpha$. It seems difficult to overcome the 
threshold $s\leq q$ giving at once expressions for $L(\chi_{t_1}\cdots\chi_{t_s},\alpha)$ with the
effectiveness of Perkins' results.

In the next Theorem, we extend the previous results beyond the mentioned threshold, providing at once new quantitative information.
\begin{Theorem}\label{functionalrelations}
Let $\alpha,s$ be positive integers, such that $\alpha\equiv s\pmod{q-1}$. Let $\delta$ be the smallest positive integer such that, simultaneously, $q^{\delta}-\alpha\geq 0$ and $s+\ell_q(q^{\delta}-\alpha)\geq 2.$ The formal series:
\begin{equation}\label{mainformula}
V_{\alpha,s}(t_1,\ldots,t_s)=\widetilde{\pi}^{-\alpha}L(\chi_{t_1}\cdots\chi_{t_s},\alpha)\omega(t_1)\cdots\omega(t_s)\prod_{i=1}^{s}\prod_{j=0}^{\delta-1}\left(1-\frac{t_i}{\theta^{q^j}}\right)\in K_\infty[[t_1,\ldots,t_s]]\end{equation} is in fact a symmetric polynomial of $K[t_1,\ldots,t_s]$ of total degree $\delta(\alpha,s)$ 
such that
$$\delta(\alpha,s)\leq s\left(\frac{s+\ell_q(q^{\delta}-\alpha)}{q-1} \right)-s.$$
\end{Theorem} This statement holds if $\alpha=q^m$ and $s\geq 2$ (so that $\delta=m$) assuming that empty products are equal to one by convention.
In this case, since $s\equiv \alpha\pmod{q-1}$, we have $s+\ell_q(q^\delta-\alpha)\equiv0\pmod{q-1}$ so that in fact, $s\geq \max\{2,q-1\}$. The reader may have noticed that the choice $\alpha=q^m$ and $s=1$ is not allowed in Theorem
\ref{functionalrelations}. However,
as mentioned above, the computation of $V_{q^m,1}$ is completely settled in \cite{Pe}. This discrimination of the 
case $\alpha=q^m,s=1$ should not be surprising neither; similarly, the Goss zeta function associated to $A$ 
has value $1$ at zero, but vanishes at all negative integers divisible by $q-1$.

In Section \ref{congruences}, we will be more specifically concerned with Bernoulli-Carlitz numbers.
A careful investigation of the polynomials $V_{1,s}$ and an application of the digit principle to the
function $\omega$ will allow us to show that, for $s\geq 2$ congruent to one modulo $q-1$,
$$\mathbb{B}_s=\Pi(s)^{-1}V_{1,s}(\theta,\ldots,\theta)$$ is a polynomial in $\FF_q[\theta]$ (Proposition \ref{theorem3}) 
(\footnote{Note that $\mathbb{B}_1$ is not well defined}). 
We shall then show the next Theorem, 
which highlights the interest of these polynomials in $\theta$.
\begin{Theorem}\label{bernoulli}
Let $s\geq 2, $ $s\equiv 1 \pmod{q-1}.$ Let us consider the expansion $s=\sum_{i=0}^rs_iq^i$ of $s$ in base $q$.  
Let $d$ be an integer such that $q^d>s$ and let $\pp$ be a prime of degree $d$. Then:
$$\mathbb{B}_s \equiv \frac{(-1)^sBC_{q^d-s}\prod_{i=0}^rl_{d-i-1}^{s_iq^i}}{\Pi(q^d-s)}\pmod{\pp}.$$
\end{Theorem}
In this result, $l_d$ denotes the polynomial $(-1)^d\prod_{i=1}^d(\theta^{q^i}-\theta)$; we observe that 
the latter polynomial is invertible modulo $\pp$ just as $\Pi(q^d-s)$.
The non-vanishing of $\mathbb{B}_s$ for fixed $s$ signifies the existence of an explicit constant $c>0$, depending on $s$ and $q$, such that
for all $d\geq c$,
\begin{equation}\label{conjcongr}
BC_{q^d-s}\not\equiv 0\pmod{\pp},\quad \text{ for all $\pp$ such that\;}\deg\pp=d.\end{equation}
However, the non-vanishing of $\mathbb{B}_s$ is also equivalent to the fact that the function $$L(\chi_{t_1}\cdots\chi_{t_s},1)\prod_{i=1}^s(t_i-\theta)^{-1},$$
entire of $s$ variables (as we will see), is a unit when identified to 
an element of $\CC_\infty[[t_1-\theta,\ldots,t_s-\theta]]$; we presently do not know how to prove this property for all $s$. 
Therefore, the property (\ref{conjcongr}) is linked with the following 
conjecture of nature analogue of classical results on the simplicity of the zeroes of Goss zeta functions and $L$-series, which should be,
we believe, true.
\begin{Conjecture}\label{conj1}
Let $s\geq 2$ be congruent to one modulo $q-1$. Then, locally at $t_1=\cdots=t_s=\theta$, 
the divisor of the zeroes of the function $L(\chi_{t_1}\cdots\chi_{t_s},1)$ is equal to the set of zeroes of the polynomial
$\prod_{i}(t_i-\theta).$
\end{Conjecture}
Numerical computations on Bernoulli-Carlitz fractions made by Taelman provide some evidence to support this 
hypothesis. The Conjecture follows from Perkins results \cite{Per2} in the case $s\leq q$ and $\alpha=s$. The conjecture is also verified
if $\ell_q(s)=q$ and $\alpha=1$, thanks to our Corollary \ref{explicitsq}.

\medskip

We end this introduction with a general remark about our methods. One of the features of this paper is the analysis of problems involving
several variables (especially in Section \ref{sectionthree}, but not only). Far from being a technical complication, this is crucial in our approach and may be 
difficult to avoid. Many corollaries we obtain by specialization in results in several variables that we obtain seem difficult to 
prove directly. 
\section{Algebraic values of the function of Anderson and Thakur\label{kroneckerweber}}

In this Section, we are going to pursue our investigation on the values of $\omega$ at the roots of unity and
we will prove Theorems \ref{mtre} and \ref{gaussthakur}. Before going on, we collect an amount of known facts and necessary notations as well as well known definitions
used all along this paper, for convenience of the reader.

\subsection{Preliminaries}

In this paper, we call {\em Carlitz's module} the unique $\FF_q$-algebra homomorphism $$\phi:A\rightarrow\mathbf{End}_{\FF_q-\text{lin.}}(\mathbb{G}_a(\mathbb{C}_\infty))$$
determined by $$\phi_\theta=\theta+\tau,$$ with $\tau$ the endomorphism such that $\tau( c)=c^q$ for all $c\in\CC_\infty$. We also recall that 
the $\FF_q$-algebra $$\mathbf{End}_{\FF_q-\text{lin.}}(\mathbb{G}_a(\CC_\infty))$$ can be identified with the skew polynomial ring $\CC_\infty[\tau]$ whose elements 
are finite sums $\sum_{i\geq 0}c_i\tau^i$ with the $c_i$'s in $\CC_\infty$, submitted to the usual product rule. We shall write $\phi_a$ for the evaluation of 
$\phi$ at $a\in A$. 

For $a\in A\setminus\{0\}$, we set $$\lambda_a=\exp\left(\frac{\widetilde{\pi}}{a}\right).$$ This is a generator
of the kernel $$\Lambda_a=\{\phi_b(\lambda_a);b\in A\}$$ of $\phi_a$ in $\CC_\infty$, an $A$-module isomorphic to $A/aA$. 

In all the following, a monic irreducible element in $A$ will be called a {\em prime}.
Let $\pp $ be a prime of $A=\FF_q[\theta]$ of degree $d$.
We denote by $K_\pp=K(\Lambda_\pp)=K(\lambda_\pp)$ the $\pp$-th cyclotomic function field extension of $K$ in $\CC_\infty$. 
We refer the reader to \cite[Chapter 12]{ROS} for the basic properties of cyclotomic function fields.
We recall here that the integral closure 
$\mathcal{O}_{K_\pp}$ of $A$ in $K_\pp$ equals the ring $A[\lambda_\pp]$.

The extension $K_\pp/K$ is cyclic of degree $q^d-1$, ramified in $\pp$ and $\theta^{-1}$. It is in fact totally ramified in $\pp$ 
and the decomposition group at $\theta^{-1}$ is isomorphic to the inertia group, therefore isomorphic to $\FF_q^\times$.
We denote by $\Delta_{\pp }$ the Galois group $\mathbf{Gal}(K(\Lambda_{\pp })/K)$. There is an unique isomorphism (Artin symbol, \cite[Proposition 7.5.4]{Go2})
$$\sigma:(A/\pp A)^\times\rightarrow\Delta_{\pp },\quad \sigma:a\mapsto\sigma_a,$$ such that
$$\sigma_a(\lambda_\pp)=\phi_a(\lambda_\pp).$$

Let $\zeta_1,\ldots,\zeta_d$ be the roots of the polynomial $\pp$ in $\FF_{q^d}$. We denote by $\FF_\pp$ the field 
$$\FF_q(\zeta_1,\ldots,\zeta_d).$$ Once chosen a root $\zeta\in\{\zeta_1,\ldots,\zeta_d\}$, the Teichm\"uller character (see \cite[Section 8.11]{Go2}) $\omega_\pp$ induces an unique group homomorphism 
$$\vartheta_\pp:\Delta_\pp\rightarrow\FF_\pp^\times,$$ defined in the following way: if $\delta=\sigma_a\in \Delta_\pp$ for some $a\in A$,
then $$\vartheta_\pp(\delta)=a(\zeta)=\chi_{\zeta}(a).$$ We will refer to this homomorphism as to the Techm\"uller character allowing an abuse 
of language (indeed, it is customary, in particular, that Teichm\"uller characters take values in Witt rings). 

\subsubsection{Gauss-Thakur sums\label{gaussthakursubs}}

For any finite abelian group $G$, we shall write $\widehat{G}$ for the group $\mathbf{Hom}(G,(\FF^{\text{alg}})^\times)$.
In particular, 
$\vartheta_\pp\in\widehat{\Delta}_\pp$. 
For the background on Gauss-Thakur sums we refer to \cite[Section 9.8]{Go2}. In our approach, however,
we find it natural to associate Gauss-Thakur sums to elements of $\widehat{\Delta}_\pp$ (compare with loc. cit. Definition 9.8.1).
\begin{Definition}{\em With $\pp$, $d$, $\vartheta_\pp$ as above,
the {\em basic Gauss-Thakur sum} $g(\vartheta_\pp^{q^j})$ associated  to this data is the element 
of $K^{\text{ab}}$ defined by:
$$g(\vartheta_\pp^{q^j})=\sum_{\delta\in\Delta_\pp }\vartheta_\pp (\delta^{-1})^{q^j}\delta(\lambda_\pp )\in\FF_\pp[\lambda_\pp].$$}\end{Definition}
The same sum is denoted by $g_j$ in \cite{Go2,Tha1}. The basic Gauss-Thakur sums are be used to define {\em general
Gauss-Thakur sums} associated to arbitrary elements of $\widehat{\Delta}_\pp$. 
The group $\widehat{\Delta}_\pp$ being isomorphic to $\Delta_\pp$ it is cyclic; it is in fact generated by $\vartheta_\pp$.
Let $\chi$ be an element of $\widehat{\Delta}_\pp$.
There exists an unique integer $i$ with $0<i<q^d$, such that
$\chi=\vartheta_\pp^i$. Let us expand $i$ in base $q$, that is, let us write $i=i_0+i_1q+\cdots+i_{d-1}q^{d-1}$ with $i_j\in\{0,\ldots,d-1\}$.
Then, $\chi=\prod_{j=0}^{d-1}(\vartheta_\pp^{q^j})^{i_j}$. 
\begin{Definition}{\em The {\em general Gauss-Thakur sum} $g(\chi)$ associated to $\chi\in\widehat{\Delta}_\pp$ as above,
is defined by:
$$g(\chi)=\prod_{j=0}^{d-1}g(\vartheta_\pp^{q^j})^{i_j}.$$}\end{Definition}

More generally, let us now consider a non-constant monic polynomial $\aa\in A$.
 We denote by $\widehat{\Delta}_\aa$ the dual character group $\mathbf{Hom}(\Delta_\aa, (\mathbb{F}_q^{\text{alg}})^\times).$ If $\chi$ is in $\widehat{\Delta}_\aa,$ we set:
$\mathbb F_q(\chi)= \mathbb F_q(\chi (\delta);\delta \in \Delta_\aa)\subset \mathbb F_q^{\text{alg}}.$
We also write $$\mathbb F_\aa=\mathbb F_q (\chi;\chi \in \widehat{\Delta}_\aa)$$ and we recall that
$\mathbf{Gal}( K_\aa(\mathbb F_\aa)/K( \mathbb F_\aa))\simeq \Delta_\aa.$
We observe that $\widehat{\Delta}_\aa$ is isomorphic to $\Delta_\aa$ if and only if $\aa$ is squarefree.
If $\aa=\pp_1\cdots\pp_n$ with $\pp_1,\ldots,\pp_n$ distinct primes, then 
$$\widehat{\Delta}_\aa\simeq \prod_{i=1}^n \widehat{\Delta}_{\pp_i}.$$
Let us then assume that $\aa$ is non-constant and squarefree. We want to extend the definition of the Gauss-Thakur sums to 
characters in $\widehat{\Delta}_\aa$. 
For $\chi \in \widehat{\Delta}_\aa,$ $\chi$ not equal to the trivial character $\chi_0,$ there exist $r$ distinct primes $\pp_1, \cdots , \pp_r$ 
and characters $\chi_1,\ldots,\chi_r$ with $\chi_j\in\widehat{\Delta}_{\pp_j}$ for all $j$, with
$$\chi=\chi_1\cdots\chi_r.$$
\begin{Definition}{\em The {\em Gauss-Thakur sum} associated to $\chi$ is the product:
$$g(\chi)=g(\chi_1)\cdots g(\chi_r).$$ The polynomial
$\mathfrak{f}_{\chi}= \pp_1\cdots \pp_r$ is called the {\em conductor} of $\chi$; it is a divisor of $\aa.$ 
The degree of $\mathfrak{f}_\chi$ will be denoted by $d_\chi$. If $\aa$ itself is a prime $\pp$ of degree $d$, then $\mathfrak{f}_{\chi}=\pp$ and $d_\chi=d$.}
\end{Definition}

The following result collects the basic properties of the sums $g(\chi)$ that we need in the sequel, and can be easily deduced from Thakur's results in \cite[Theorems I and II]{Tha1}.
\begin{Proposition}\label{Thakur} Let $\aa\in A$ be monic, squarefree of degree $d$. 
The following properties hold.
\begin{enumerate}
\item If $\chi=\chi_0$ is the trivial character, then $g(\chi)=1$.
\item For all $\delta \in \Delta_\aa,$ we have $\delta (g (\chi ))=\chi (\delta) g (\chi).$
\item If $\chi \not = \chi_0, $ then $g(\chi)g(\chi^{-1})=(-1)^{d_{\chi}}\mathfrak{f}_{\chi}.$
\end{enumerate}
\end{Proposition}
By the normal basis theorem, $K_\aa$ is a free $K[\Delta_\aa]$-module of rank one.
Gauss-Thakur's sums allow to determine explicitly a generator of this module:
\begin{Lemma}
\label{normalbasis} Let us write $\eta_\aa=\sum_{\chi\in \widehat{\Delta}_\aa}g(\chi) \in K_\aa.$ Then :
$$K_\aa=K[\Delta_\aa]\eta_\aa,$$
and 
$$ A_\aa=A[\Delta_\aa]\eta_\aa,$$
where $ A_\aa$ is the integral closure of $A$ in $K_\aa.$

Moreover, let $\chi$ be in $\widehat{\Delta}_\aa$. Then, the following identity holds:
\begin{equation}\label{eqkffaa}
K(\FF_{\aa})g(\chi)=\{ x\in K_\aa(\mathbb F_\aa)\text{ such that for all }\delta \in  \Delta_\aa, \delta (x)=\chi(\delta) x\}.
\end{equation}
\end{Lemma}
\noindent\emph{Proof.} Let us expand $\aa$ in product $\pp_1\cdots\pp_n$ of distinct primes $\pp_i$. To show that
$ A_\aa=A[\Delta_\aa]\eta_\aa$ (this yields the identity $K_\aa=K[\Delta_\aa]\eta_\aa$) one sees that
$$ A_\aa\simeq  A_{\pp_1} \otimes _A\cdots  \otimes_A  A_{\pp_n},$$ because the discriminants of the extensions
$ A_{\pp_i}/A$ are pairwise relatively prime and the fields $K_{\pp_i}$ are pairwise linearly disjoint (see \cite{FrTa}).
One then uses \cite[Th\'eor\`eme 2.5]{ANG} to conclude with the second identity. 

We now prove the identity (\ref{eqkffaa}).
We recall that if we set, for $\chi\in\widehat{\Delta}_\aa$,
$$e_\chi=\frac{1}{|\Delta_\aa|}\sum_{\delta\in\Delta_\aa}\chi(\delta)\delta^{-1}\in\FF_q(\chi)[\Delta_\aa]$$
(well defined because $p$, the rational prime dividing $q$, does not divide $|\Delta_\aa|$), then the following identities hold:
\begin{itemize}
\item $e_\chi e_\psi=\delta_{\chi,\psi}e_\chi$ (where $\delta_{\chi,\psi}$ denotes Kronecker symbol),
\item for all $\delta\in\Delta_\aa$, $\delta e_\chi=\chi(\delta)e_\chi$,
\item $\sum_{\chi\in\widehat{\Delta}_\aa}e_\chi=1$.
\end{itemize}
This yields $e_\chi\eta_\aa=g(\chi)$. Now, by $K(\FF_\aa)=K_\aa(\FF_\aa)[\Delta_\aa]\eta_\aa$, we get
$e_\chi K_\aa(\FF_\aa)=K(\FF_\aa)g(\chi)$. The second part of the Lemma then follows by observing that
if $M$ is an $\FF_\aa[\Delta_\aa]$-module, then
$$e_\chi M=\{m\in M\text{ such that for all }\delta\in\Delta_\aa,\delta m=\chi(\delta)m\}.$$
\CVD

\subsubsection{The function of Anderson and Thakur}

For the basic properties of this 
function, introduced in \cite{AT} (see also \cite{Bourbaki}), we also suggest to read \cite[Section 3.1]{gamma}. 
We recall from \cite[Corollaries 5, 10]{Pe} that $t\in\CC_\infty\setminus\{\theta^{q^k};k\geq 0\}$ and $\omega(t)$
are simultaneously algebraic if and only if $t=\zeta\in\FF_q^{\text{alg}}$. 
Moreover, we have the following Lemma, where we adopt the convention that $\FF_{q^0}:=\emptyset$.
\begin{Lemma}\label{propertiesomega} Let $d\geq 1$ be an integer.
For $\zeta\in\FF_{q^d}\setminus\FF_{q^{d-1}}$, $\omega(\zeta)$ belongs to the set
$$\FF_{q^d}^\times\rho_\zeta,$$ where 
$\rho_\zeta$ is a distinguished root of the polynomial 
\begin{equation}\label{thepolynomial}
X^{q^d-1}-(\zeta-\theta^{q^{d-1}})\cdots(\zeta-\theta)\in A[\zeta][X].\end{equation}
Moreover, regardless of the choice of $\zeta$,
$$v_\infty(\omega(\zeta))=-\frac{1}{q-1}.$$
\end{Lemma}
\noindent\emph{Proof.} This is a simple consequence of \cite[Corollary 5]{Pe}.\CVD
The function of Anderson and Thakur can also be defined, alternatively, by the series expansion:
\begin{equation}\label{scarlitz}
\omega(t):=\sum_{i=0}^\infty \lambda_{\theta^{i+1}}t^i=\sum_{n=0}^\infty\frac{\widetilde{\pi}^{q^n}}{d_n(\theta^{q^n}-t)}
\in (-\theta)^{1/(q-1)}K_\infty[[t]]
\end{equation} converging for $|t|<q$. 

The {\em Tate algebra} $\TT_t$ is the $\CC_\infty$-algebra 
whose elements are the series 
$\sum_{i\geq 0}c_it^i\in\CC_\infty[[t]]$ converging in the bordered unit disk $$\overline{D}(0,1)=\{t\in\CC_\infty,|t|\leq1\}.$$
Here, $|\cdot|$ denotes an absolute value associated to the valuation $v_\infty$. In all the following, for clarity, 
we normalize it by setting $|\theta|=q$.

We recall from \cite[Section 4]{Pe} that the Tate algebra $\TT_t$ is endowed with the norm 
$\|\cdot\|$ defined as follows: if $f=\sum_{i\geq 0}c_it^i\in\TT_t$ with $c_i\in\CC_\infty$ ($i\geq 0$), then
$\|f\|=\sup_{i\geq 0}|c_i|=\max_{i\geq 0}|c_i|$. Endowed with this norm, $\TT_t$ becomes a $\CC_\infty$-Banach algebra.
For $r>0$ a real number, we denote by $\mathbb{D}_r$ the $\FF_q[t]$-submodule of $\TT_t$
whose elements are the series $f$ such that $\|f\|<r$. The operator $\tau$ extends in an unique way to a $\FF_q[t]$-automorphism
of $\TT_t$ so we have at once all the $\FF_q[t]$-endomorphisms $\phi_{\mathfrak{a}}-\chi_t(\mathfrak{a})$.

In all the following, we denote by $K[t][\tau]$ and $K[[\tau]]$ respectively the skew polynomial rings in powers of $\tau$ with coefficients in $K[t]$
and the skew entire series rings in powers of $\tau$ with coefficients in $K$, endowed with the product rule induced by the identity $\tau t=t\tau$.
The evaluation operator (see \cite[Section 4]{Pe}) $$E_\mathfrak{e}:\CC_\infty\rightarrow\CC_\infty,$$
where
$$\mathfrak{e}=\sum_{i\geq 0}\frac{\tau^i}{d_i}\in K[[\tau]]$$ is the series associated to Carlitz's exponential function introduced in 
\cite[Section 4]{Pe},
extends to an $\FF_q[t]$-endomorphism of $\TT_t$. The above formula (\ref{scarlitz})
can be rewritten in a compact form as
$$\omega(t)=E_{\mathfrak{e}}\left(-\frac{\widetilde{\pi}}{t-\theta}\right).$$ We mention that in \cite{gamma}, following a suggestion of D. Goss, an
analogy between the above formula and the definition of the gamma function as a Mellin transform of the function $e^{-z}$ was discussed.

From (\ref{scarlitz}), one deduces easily that $\omega$ belongs to $\TT_t$ (see Subsection \ref{omegaaj} below).
It is also easy to show that $\omega$ is a generator of the free $\FF_q[t]$-module of rank one, kernel of 
the operator $$\phi_\theta-\chi_t(\theta)=\tau+\theta-t\in K[t][\tau],$$
so that
$$\tau\omega(t)=(t-\theta)\omega(t)$$
 (see \cite[Proposition 3.3.6]{Pa}).
This implies that $\omega$ is a generator of the intersection of the kernels in $\TT_t$ of the various
operators $$\phi_\aa-\chi_t(\aa ),$$ with $\aa$ monic as above (see \cite[Lemma 29]{Pe}). We are going 
to provide, in Proposition \ref{basis}, a complete description of the kernels of each one of these operators.

\subsection{The functions $\omega_{\aa,j}$\label{omegaaj}}

The following Lemma holds.
\begin{Lemma}\label{lemmaevalexp} The kernel of the evaluation operator $E_\mathfrak{e}$ is
the submodule $\widetilde{\pi}A[t]$ of $\TT_t$. Its restriction to $\mathbb{D}_{|\widetilde{\pi}|}$ is an isometry.
\end{Lemma}
\noindent\emph{Proof.} The kernel clearly contains $\widetilde{\pi}A[t]$.
Let $m$ be an element of $\text{ker}(E_\mathfrak{e})$. Then, $m=\sum_{i\geq 0}c_it^i$ with $c_i\in\CC_\infty$ and 
$\exp(c_i)=0$, so that $c_i\in\widetilde{\pi}A$ for all $i$. But then, 
$$\text{ker}E_\mathfrak{e}\subset \widetilde{\pi}A[[t]]\cap\TT_t=\widetilde{\pi}A[t].$$
That this endomorphism is an isometry on $\mathbb{D}_{|\widetilde{\pi}|}$ was implicitly observed in \cite{Pe}. 
This relies on the fact that $\exp$ induces an isometry on the disk $\{z\in \CC_\infty;|z|<|\widetilde{\pi}|\}$ and the simple verification is left to the reader.\CVD

In order to define the functions $\omega_{\aa,j}$ we will compute, following \cite[Section 4]{Pe}, 
the image of $E_\mathfrak{e}$ at various rational functions of $\TT_t$, and for this, we will need the 
next Lemma.
\begin{Lemma}\label{rationalfunctions}
Let $\aa$ be a non-constant monic polynomial of $A$. Then, $1/(\aa-\chi_t(\aa))\in\TT_t$.
\end{Lemma}
\noindent\emph{Proof.} It suffices to show that the roots $\xi\in\CC_\infty$ of the polynomial $\aa-\chi_t(\aa)\in A[t]$ are all such that $|\xi|>1$.
But this is obvious, since we have $|\chi_\xi(\aa)|>1$ if and only if $|\xi|>1$.\CVD

We now fix a monic polynomial $\aa$ of degree $d>0$. 
Lemma \ref{rationalfunctions} implies that, for $j=0,\ldots,d-1$, the series
 $$\omega_{\aa,j}(t)=E_{\mathfrak{e}}\left(\frac{\theta^j\widetilde{\pi}}{\aa -\chi_t(\aa)}\right)$$ are well defined elements of $\TT_t$.
 When the reference to the polynomial $\aa$ is clear, we will write $\omega_j$ instead of $\omega_{\aa,j}$. In particular, $\omega_{\theta,0}=\omega$.
 By Lemma \ref{lemmaevalexp}, we have
 \begin{equation}\label{norms}
 \|\omega_j\|=\left\|\frac{\theta^j\widetilde{\pi}}{\aa-\chi_t(\aa)}\right\|=|\widetilde{\pi}\theta^ja^{-1}|=|\phi_{\theta^j}(\lambda_\aa)|=q^{\frac{q}{q-1}+j-d},\quad j=0,\ldots,d-1.
 \end{equation}
 
To study the elements $\omega_j$ as rigid analytic functions, it may be convenient to observe that the function $\chi_t(\aa)$ (with $\aa$ as above)
induces a rigid analytic endomorphism of the bordered unit disk $\overline{D}(0,1)=\{t\in\CC_\infty;|t|\leq 1\}$. Hence, 
right composition in series of powers of a new variable $x$ by setting $x=\chi_t(\aa)$ induces a map $\TT_x\rightarrow\TT_t$.

Assuming that $\tau$ is linearly extended to $\TT_x$ by the rule $\tau(x)=x$, let us now consider the series
$$\omega_{j}^*(x)=\omega_{\aa,j}^*(x)=E_{\mathfrak{e}}\left(\frac{\widetilde{\pi}}{\aa-x}\right)=\sum_{n\geq 0}\exp\left(\frac{\widetilde{\pi}}{\aa^{n+1}}\right)x^{n}
\in\TT_x,\quad j=0,\ldots,d-1.$$ Again, 
$\omega^*_{\theta,x}(x)=\omega(x)$. Furthermore, we notice that, for all $j=0,\ldots,d-1$, $\omega_j^*$ has no zeroes on the disk $\overline{D}(0,1)$.
Then, 
\begin{equation}\label{idcomposition}
\omega_j(t)=\omega^*_j\circ\chi_t(\aa),\quad j=0,\ldots,d-1,
\end{equation}
and we see that these functions have no zeroes in the disk $\overline{D}(0,1)$.
 
\subsubsection{Kernel of the operators $\phi_{\aa}-\chi_t(\aa)$}
 
We shall prove the next Proposition.

 \begin{Proposition}\label{basis} The kernel in $\TT_t$ of the operator $\phi_\aa-\chi_t(\aa)$ is the free $\FF_q[t]$-module of
rank $d$ generated by the series $\omega_{\aa,0},\ldots,\omega_{\aa,d-1}$.
\end{Proposition}

\noindent\emph{Proof.} We will write $\omega_j$ at the place of $\omega_{\aa,j}$ for simplicity. Let us consider the column matrix
$$\Omega_{\aa}(t)=\left(\begin{array}{c}\omega_0\\ \omega_1\\ \vdots\\ \omega_{d-1}\end{array}\right)\in\mathbf{Mat}_{d\times 1}(\TT_t).$$
Let us write:
$$\aa =a_0+a_1\theta+\cdots+a_{d-1}\theta^{d-1}+\theta^d\in A^+,\quad a_0,\ldots,a_{d-1}\in\FF_q.$$
By the identity
\begin{eqnarray*}
E_{\mathfrak{e}}\left(\frac{\theta^d\widetilde{\pi}}{\aa -\chi_t(\aa)}\right)&=&E_{\mathfrak{e}}\left(\frac{(\aa -\chi_t(\aa)+\chi_t(\aa)-a_0-a_1\theta-\cdots-a_{d-1}\theta^{d-1})\widetilde{\pi}}{\aa -\chi_t(\aa)}\right)\\
&=&(\chi_t(\aa)-a_0)\omega_0-a_1\omega_1-\cdots-a_{d-1}\omega_{d-1},
\end{eqnarray*}
We obtain
$$\phi_\theta\Omega_{\aa}=M_{\aa }(t)\Omega_{\aa},$$
where 
$$M_{\aa}(t)=\left(\begin{array}{ccccc}0&1&0&\cdots&0\\ 0&0&1&\cdots&0\\  
\vdots & \vdots & \vdots & & \vdots \\
0&0&0&\cdots&1\\
-a_0 & -a_1 & -a_2 & \cdots & -a_{d-1}\end{array}\right)+\left(\begin{array}{cccc}0&0&\cdots&0\\ 0&0&\cdots&0\\  
\vdots & \vdots & & \vdots \\
0&0&\cdots&0\\
\chi_t(\aa) & 0 & \cdots & 0\end{array}\right)\in\mathbf{GL}_{d}(\FF_q(t)),$$ compare with \cite[Section 5.3]{Go2}.
Moreover, $M_\aa (t)$ commutes with $\tau$, therefore, the matrix $\aa (M_\aa (t))$ represents the scalar multiplication by $\chi_t(\aa)$ and
$$\phi_\aa\Omega_\aa =\chi_t(\aa)\Omega_\aa .$$ This already shows that $\omega_0,\ldots,\omega_{d-1}$ belong
to the kernel of $\phi_\aa-\chi_t(\aa)$.

We now show that these functions are linearly independent over $\FF_q[t]$.
Let us assume by contradiction that there exist elements $\mu_0,\ldots,\mu_{d-1}\in\FF_q[t]$, not all zero, such that
\begin{equation}\label{linearindep}
\sum_{i=0}^{d-1}\mu_i(t)\omega_i(t)=0.
\end{equation} We may even assume, without loss of generality, that there exists a root $\zeta$ of $\aa$,
and an index $i\in\{0,\ldots,d-1\}$, such that 
$\mu_i(\zeta)\neq0$. By (\ref{linearindep}), we see that 
$$E_{\mathfrak{e}}\left(\frac{\widetilde{\pi}}{\aa-\chi_t(\aa)}\sum_{i=0}^{d-1}\mu_i(t)\theta^i\right)=0$$ and, by Lemma \ref{lemmaevalexp},
there exists an element $b\in A[t]$ such that 
$$\sum_{i=0}^{d-1}\mu_i(t)\theta^i=b(t)(\aa-\chi_t(\aa)).$$ Evaluating at $t=\zeta$ now yields:
$$\sum_{i=0}^{d-1}\mu_i(\zeta)\theta^i=b(\zeta)\aa.$$ The above hypothesis on the $\mu_i$'s implies that $b(\zeta)\neq0$. 
However, the degrees in $\theta$ of the left- and right-hand sides do not agree, in contradiction with our assumption, so that
$\omega_0,\ldots,\omega_{d-1}$ are linearly independent over $\FF_q[t]$.

To finish the proof of the Proposition, we still need to show that the kernel of the operator $\phi_\aa-\chi_{t}(\aa)$
is the free $\FF_q[t]$-module of rank $d$ generated by the functions $\omega_0,\ldots,\omega_{d-1}$. Now, $\phi_\aa-\chi_{t}(\aa)$
operates on the fraction field $\mathbb{L}_t$ of $\TT_t$ as well. The subfield of $\LL_t$ of elements 
which are $\tau$-invariant is equal to $\FF_q(t)$. Therefore, the kernel of $\phi_\aa-\chi_{t}(\aa)$ is a $\FF_q(t)$-vector space of 
dimension $\leq d$ (because of the Wronskian Lemma). But $\omega_0,\ldots,\omega_{d-1}$ are linearly independent
over $\FF_q[t]$, hence over $\FF_q(t)$, and belong to the kernel which then is equal to the $\FF_q(t)$-vector space
generated
by $\omega_0,\ldots,\omega_{d-1}$. 

Let $f\in\TT_t$ be such that
$\phi_\aa(f)-\chi_{t}(\aa)f=0$ and let us consider an element $g$ in $\TT_t$ such that
$E_{\mathfrak{e}}(g)=f$. By Lemma \ref{lemmaevalexp}, we have 
$(a-\chi_t(a))g\in\widetilde{\pi}A[t]$ so that 
$$g\in \widetilde{\pi} A[t]+\sum_{j=0}^{d-1}\FF_q[t]\frac{\theta^j\widetilde{\pi}}{a-\chi_t(a)}.$$
Therefore, 
$$f=E_{\mathfrak{e}}(g)\in\sum_{j=0}^{d-1}\FF_q[t]\omega_j(t).$$\CVD

\begin{Remark}{\em In analogy with 
H\"older's Theorem for the gamma function, it is not difficult to show that
the functions $$\omega,\mathcal{D}_1\omega,\ldots,\mathcal{D}_n\omega,\ldots$$ are algebraically independent over $\CC_\infty(t)$.}
\end{Remark}

\begin{Remark}{\em Let $n$ be a nonnegative integer, let us consider a prime $\mathfrak{p}\in A$ of degree $d\geq 1$. It is possible to show, with the same ideas as in the proof of \cite[Proposition 15]{gamma},
that the kernel of the operator $$\phi_{\mathfrak{p}^{n+1}}-\chi_t(\mathfrak{p}^{n+1}),$$ of dimension $nd$ where $d$ is the degree of $\pp$, 
is spanned by the entries of the matrices
$$\Omega_\pp,\mathcal{D}_1\Omega_\pp,\ldots,\mathcal{D}_n\Omega_\pp.$$}\end{Remark}
 
 \subsection{Analytic identities}
 
 In this subsection we proceed to review the main properties of the Tate algebra and higher derivatives we need and then, we 
 describe two families of analytic identities that will be of crucial use in the proof of Theorem \ref{mtre}.
 
 \subsubsection{Tate algebras and higher derivatives\label{sectchainrule}}

Let $t_1,\ldots,t_s$ be independent variables. Extending the observations of the beginning of
Section \ref{omegaaj}, we consider now the {\em Tate algebra} in $s$ variables $\TT_{t_1,\ldots,t_s}$, that is,
the $\CC_\infty$-subalgebra of $\CC_\infty[[t_1,\ldots,t_s]]$ whose elements $f$ are formal series
\begin{equation}\label{specimen}
f=\sum_{i_1,\ldots,i_s\in\ZZ_{\geq 0}}f_{i_1,\ldots,i_s}t_1^{i_1}\cdots t_s^{i_s},\quad f_{i_1,\ldots,i_s}\in\CC_\infty
\end{equation}
converging in the bordered unit polydisk $\overline{D}(0,1)^s$.
The algebra $\TT_{t_1,\ldots,t_s}$ is endowed with the norm $\|\cdot\|$ generalizing the norm used in Section \ref{omegaaj} and defined as follows.
Let $f\in \TT_{t_1,\ldots,t_s}$ as in (\ref{specimen}). Then,
$$\|f\|:=\sup_{i_1,\ldots,i_s}|f_{i_1,\ldots,i_s}|=\max_{i_1,\ldots,i_s}|f_{i_1,\ldots,i_s}|.$$
We recall that with this norm, $\TT_{t_1,\ldots,t_s}$ is a Banach $\CC_\infty$-algebra. Furthermore, in this setting,
we can extend $\tau$ to a $\FF_q[t_1,\ldots,t_s]$-automorphism and we have the 
$\FF_q[t_1,\ldots,t_s]$-endomorphism $E_{\mathfrak{e}}$. 

It is helpful to also notice that if 
$f_0,f_1,\ldots,f_s$ are elements of $\TT_{t_1,\ldots,t_s}$ with $\|f_i\|\leq 1$ for $i=1,\ldots,s$, then the composition of 
functions
$f_0(f_1,\ldots,f_s)$ also defines an element of $\TT_{t_1,\ldots,t_s}$.
Let $\TT_{t_1,\ldots,t_d}^\text{sym}$ be the sub-algebra of $\TT_{t_1,\ldots,t_d}$ of functions which are symmetric 
in the variables $t_1,\ldots,t_d$. The automorphism $\tau$ induces a $\FF_q[t_1,\ldots,t_d]^{\text{sym}}$-automorphism
of $\TT_{t_1,\ldots,t_d}^\text{sym}$, where we have denoted by $\FF_q[t_1,\ldots,t_d]^{\text{sym}}$ the subring of 
polynomials which are symmetric in $t_1,\ldots,t_d$.

We will need the following Lemma.
\begin{Lemma}\label{symmetric1}
Let $\pp$ be a prime of $A$ of degree $d$, let $\zeta_1,\ldots,\zeta_d$ be its roots in $\FF_q^{\text{alg}}$, let us consider 
an element $f\in \TT_{t_1,\ldots,t_d}^\text{sym}$. Then,
$$(\tau f)(\zeta_1,\ldots,\zeta_d)=f(\zeta_1,\ldots,\zeta_d)^q.$$
\end{Lemma}
\noindent\emph{Proof.} By virtue of \cite[Th\'eor\`eme 2, A IV.64]{Bou}, 
the series of $f$ can be expanded as a series in the elementary symmetric polynomials
$$\alpha_n(t_1,\ldots,t_d)=\sum_{I\subset\{1,\ldots,d\}}\left(\prod_{i\in I}t_i\right),\quad n=1,\ldots,d$$
(the sum running over the subsets $I$ of $\{1,\ldots,d\}$ of size $n$), we only need to verify the Lemma 
for $f=\alpha_n$, $n=0,\ldots,d$ (indeed, $\CC_\infty[t_1,\ldots,t_d]^{\text{sym}}$ is dense for the sup-norm $\|\cdot\|$
in $\TT_{t_1,\ldots,t_d}^{\text{sym}}$). But in this case,
$$(\tau \alpha_n)(\zeta_1,\ldots,\zeta_d)=\alpha_n(\zeta_1,\ldots,\zeta_d)=\alpha_n(\zeta_1,\ldots,\zeta_d)^q,$$ because 
$\zeta_1,\ldots,\zeta_d$ are conjugate.\CVD

\noindent\emph{Higher derivatives.}
We will occasionally need to compute the higher derivatives of $\omega$ and other allied functions.
For the background on higher derivatives (also called {\em hyperderivatives}), we refer to the recent work of Jeong \cite{Je}
noticing that the specific tools we are interested in are also contained in Teichm\"uller's paper \cite{Tei}.

The $\CC_\infty$-linear higher derivative $(\mathcal{D}_{t,n})_{n\geq 0}$ (also denoted by $(\mathcal{D}_{n})_{n\geq 0}$
in this text) defined by 
$$\mathcal{D}_{t,m}(t^n)=\binom{n}{m}t^{n-m},\quad m,n\geq0$$ induces $\CC_\infty$-linear endomorphisms of 
$\TT_t$. As an example of computation, we have
\begin{equation}\label{eqdn}
\mathcal{D}_{t,n}(\theta-t)^{-1}=(\theta-t)^{-n-1},\quad n\geq 0.\end{equation}

This higher derivative can be defined over more general $\CC_\infty$-algebras
of rigid analytic functions of the variable $t$ as well.
It also satisfies the chain rule, see \cite[Section 2.2]{Je} and 
\cite[Equation (6)]{Tei}. In other words, for all $n\geq 1$, there exist 
polynomials 
$$F_{n,i}(X_1,\ldots,X_{n+1-i})\in\FF_p[X_1,\ldots,X_{n+1-i}],\quad i=1,\ldots,n$$ (where $p$ is the prime dividing $q$) with the following property.
For $f,g$ rigid analytic functions with $g$ defined over a non-empty open subset $\mathcal{O}\subset\CC_\infty$
and $\CC_\infty$-valued in such a way that $f\circ g$ is a well defined rigid analytic function $\mathcal{O}\rightarrow\CC_\infty$,
\begin{equation}\label{chainrule}
\mathcal{D}_{t,n}(f\circ g)=\sum_{i=1}^nF_{n,i}(\mathcal{D}_{t,1}g,\ldots,\mathcal{D}_{t,n+1-i}g)(\mathcal{D}_{t,i}f)\circ g.\end{equation} Moreover,
one easily sees that
$$F_{n,n}=X_1^n,\quad F_{n,1}=X_n,\quad n\geq 1.$$
This property holds in particular in $\TT_t$ for $f,g\in\TT_t$ and $\|g\|\leq 1$, when $f\circ g\in\TT_t$. 
Clearly, for all $n$, $\mathcal{D}_n$ commutes with the operator $E_\mathfrak{e}$. We will 
often write $\mathcal{D}_{n}$ instead of $\mathcal{D}_{t,n}$ to simplify our notations. 
More generally, for all $i=1,\ldots,s$, $(\mathcal{D}_{t_i,n})_{n\geq 0}$ is a higher derivative of $\TT_{t_1,\ldots,t_s}$.

\subsubsection{A first family of analytic identities} 
 
 The next result we need is Proposition \ref{anidentity} below.
  For $i=0,\ldots,d-1$, we set:
$$\aa_i=a_{i+1}+a_{i+2}\theta+\cdots+a_{d-1}\theta^{d-i-2}+\theta^{d-i-1},$$
so that, in particular, $\aa_1=a_1+a_2\theta+\cdots+a_{d-1}\theta^{d-2}+\theta^{d-1}$, $\aa_{d-2}=a_{d-1}+\theta$,
$\aa_{d-1}=1$. We also set $\aa_{-1}:=\aa $ for completeness.
 \begin{Proposition}[First family of analytic identities]\label{anidentity}
The following identity holds in $\TT_t$:
\begin{equation}\label{identity}\omega(t)=\sum_{i=0}^{d-1}\chi_t(\aa_{i})\omega_i(t).\end{equation}
Moreover, for all $n\geq 1$, there exists an element $\Omega_n$ of the submodule of $\TT_t$
$$\sum_{i=0}^{n-1}\sum_{j=0}^{d-1}\FF_q[t](\mathcal{D}_i\omega_j)$$
such that
\begin{equation}\label{anidentityn}
(\mathcal{D}_n\omega)(t)=\sum_{j=0}^{d-1}\chi_t(\aa_j)(\mathcal{D}_n\omega_{j})(t)+\Omega_{n}.
\end{equation}
\end{Proposition}
\noindent\emph{Proof.} 
In $\FF_q[t,\theta]$, we have the elementary identity:
\begin{equation}\label{for0}
\frac{\aa-\chi_t(\aa)}{\theta-t}=\sum_{j=0}^{d-1}\chi_{t}(\aa_j)\theta^{j}. 
\end{equation}
Multiplying both sides of it by $\widetilde{\pi}$ and dividing by $\aa-\chi_t(\aa)$, we obtain the identity 
$$\frac{\widetilde{\pi}}{\theta-t}=\sum_{j=0}^{d-1}\chi_{t}(\aa_j)\frac{\theta^{j}\widetilde{\pi}}{\aa-\chi_t(\aa)},
$$
 which holds in $\TT_t$ by Lemma \ref{rationalfunctions}. The first part of the proposition follows after evaluation of both sides by $E_{\mathfrak{e}}$.
 
 The second part of the proposition is a direct consequence of Leibniz formula for $(\mathcal{D}_n)_{n\geq 0}$, but we give all the details
 of the intermediate computations for convenience of the reader.
 Let us consider again (\ref{for0}) and apply $\mathcal{D}_{n}$ on both left- and right-hand sides. By (\ref{eqdn}) and the chain rule (Subsection \ref{sectchainrule}), we verify that
 \begin{equation}\label{dnfraction}
 \mathcal{D}_{n}\left(\frac{1}{\aa-\chi_t(\aa)}\right)=\frac{\chi_t(\aa')^n}{(\aa-\chi_t(\aa))^{n+1}}+\Xi_n,\end{equation}
 where, by Leibniz rule, $\Xi_n$ is an element of the submodule of $\TT_t$:
 \begin{equation}\label{submodule}
 \sum_{j=0}^{n-1}\FF_q[t]\mathcal{D}_j\left(\frac{1}{a-\chi_t(a)}\right).\end{equation}
 Therefore, by  (\ref{for0}),
 $$\frac{1}{(\theta-t)^{n+1}}=\sum_{j=0}^{d-1}\theta^j\left(\frac{\chi_t(\aa_j)\chi_t(\aa')^n}{(\aa-\chi_t(\aa))^{n+1}}+\Upsilon_{j,n}\right),$$ where $\Upsilon_{j,n}$ again are elements of the submodule (\ref{submodule}). These identities hold in $\TT_t$
 and since $\mathcal{D}_n$ commutes with $E_\mathfrak{e}$, we get, multiplying by $\widetilde{\pi}$ and applying $E_\mathfrak{e}$:
 \begin{eqnarray*}
 (\mathcal{D}_n\omega)(t)&=&\sum_{j=0}^{d-1}\chi_t(\aa_j)\chi_t(\aa')^nE_\mathfrak{e}\left(\frac{\widetilde{\pi}\theta^j}{(\aa-\chi_t(\aa))^{n+1}}\right)+\Omega_{j,n}\\
 &=&\sum_{j=0}^{d-1}\chi_t(\aa_j)\chi_t(\aa')^n(\mathcal{D}_n\omega_{\aa,j}^*)(\chi_t(\aa))+\Omega_{j,n}\\
 &=&\sum_{j=0}^{d-1}\chi_t(\aa_j)(\mathcal{D}_n\omega_{\aa,j})(t)+\Omega_{j,n}\\
 \end{eqnarray*}  
 where $\Omega_{j,n}$ is an element of 
 $$\FF_q[t]\omega_{\aa,j}+\cdots+\FF_q[t](\mathcal{D}_{n-1}\omega_{\aa,j}).$$
 \CVD
  
\subsubsection{A second family of analytic identities}

Our second family of identities holds in $\TT_{t_1,\ldots,t_d}^\text{sym}$, for $d\geq 1$. 
Let us consider the 
higher derivatives $\mathcal{D}_{t_i,n}$ ($i=1,\ldots,d$)
of $$R=\mathbb{F}_q[t_1,\ldots,t_d].$$
Their {\em sum}
$$\partial=(\partial_n)_{n\geq0}$$ is the family of operators defined by:
$$\partial_n=\sum_{k_1+\cdots+k_d=n}\mathcal{D}_{t_1,k_1}\cdots\mathcal{D}_{t_d,k_d},\quad n\geq 0.$$
It is easy to verify that this also is a higher derivative (to check this, one can use, for example, the multinomial theorem). 
Their sum also induces a higher derivative of $R^{\text{sym}}(\theta)$, of its completion with respect to the 
ideal $(t_1,\ldots,t_d)$ (the ring $K[[t_1,\ldots,t_d]]^{\text{sym}}$ of symmetric formal series in powers of $t_1,\ldots,t_d$ with coefficients in $K$),
as well as of $\mathbb{T}_{t_1,\ldots,t_d}^{\text{sym}}$. We also notice that the polynomial $P=\prod_{i=1}^d(\theta-t_i)\in R^{\text{sym}}[\theta]$ is a unit
of $\mathbb{T}_{t_1,\ldots,t_d}^{\text{sym}}$. We will need the next elementary Lemma.
 \begin{Lemma}\label{fourpropoerties}
Let $N$ be a polynomial of $R^{\text{sym}}[\theta]$. For all $n\geq 0$, we have, in $\mathbb{T}_{t_1,\ldots,t_d}^{\text{sym}}$,
$$\partial_n\left(\frac{N}{P}\right)=\frac{(\mathcal{D}_{\theta,1}P)^nN}{P^{n+1}}+\Psi_n,$$
where $\Psi_n$ is an element of the module $R^{\text{sym}}[\theta]P^{-n}$. \end{Lemma}
\noindent\emph{Proof.} It is obvious that, for all $j\leq d$, $\partial_n(\theta-t_j)^{-1}=(-1)^n\mathcal{D}_{\theta,n}(\theta-t_j)^{-1}$. 
Therefore, 
$$\partial_nP^{-1}=(-1)^n\mathcal{D}_{\theta,n}P^{-1},$$
so that, by Leibniz's rule and the chain rule, 
\begin{eqnarray*}
\partial_n\left(\frac{N}{P}\right)&=&(-1)^n(\mathcal{D}_{\theta,n}P^{-1})N+\Theta_n\\
&=&\frac{(\mathcal{D}_{\theta,1}P)^nN}{P^{n+1}}+\Psi_n
\end{eqnarray*} ($\Theta_n$ is another element of $R^\text{sym}[\theta]P^{-n}$).\CVD
For $\aa\in A$ we denote by   $U_\aa$ the polynomial of $R^{\text{sym}}[\theta]$ such that:
\begin{equation}\label{primitiveformula}
\frac{U_{\aa}}{P}=\sum_{i=1}^d\frac{\aa(t_i)}{\theta-t_i}.\end{equation} The next Proposition holds:
\begin{Proposition}[Second family of analytic identities]\label{anidentitytwo}
The following identity holds in the algebra $\TT_{t_1,\ldots,t_d}^{\text{sym}}$:
\begin{equation}\label{secondfamily1}\sum_{i=1}^d\chi_{t_i}(\aa)\omega(t_i)=E_\mathfrak{e}\left(\frac{\widetilde{\pi}U_\aa}{P}\right).
\end{equation}
More generally, for all $n\geq 0$, we have the identity in $\TT_{t_1,\ldots,t_d}^{\text{sym}}$:
\begin{equation}\label{secondfamily2}
\sum_{i=1}^d\aa(t_i)(\mathcal{D}_{n}\omega)(t_i)=E_\mathfrak{e}\left(\frac{\widetilde{\pi}(\mathcal{D}_{\theta,1}P)^nU_\aa}{P^{n+1}}\right)+\Lambda_n,
\end{equation}
where $\Lambda_n$ is an element of 
$$R^{\text{sym}}E_\mathfrak{e}\left(\frac{\widetilde{\pi}\FF_q[\theta]}{P^n}\right).$$
\end{Proposition}
\noindent\emph{Proof.} The first identity follows by multiplying both sides of (\ref{primitiveformula}) by $\widetilde{\pi}$
and applying $E_\mathfrak{e}$. We now prove the second identity.
For all $n\geq 0$, we have the next identity in $\TT_{t_1,\ldots,t_d}^{\text{sym}}$, in virtue of Lemma \ref{fourpropoerties}:
\begin{equation}\label{firstintermediatestep}
\partial_{n}\left(\frac{U_\aa}{P}\right)=\frac{(\mathcal{D}_{\theta,1}P)^nU_\aa}{P^{n+1}}+\Gamma_n,
\end{equation}
with $\Gamma_n$ an element of the module $R^\text{sym}[\theta]P^{-n}$.
Hence, applying $\partial_n$ on both sides of (\ref{primitiveformula}) we get
\begin{equation}\label{nexttolast}
\sum_{i=1}^d\frac{a(t_i)}{(\theta-t_i)^{n+1}}=\frac{(\mathcal{D}_{\theta,1}P)^nU_\aa}{P^{n+1}}+\Gamma_n+\Sigma_n,
\end{equation}
where $\Sigma_n$ is an element of
$$\sum_{i=1}^d\FF_q[t_i]\frac{1}{(\theta-t_i)^n}\cap\TT_{t_1,\ldots,t_d}^{\text{sym}},$$ module which is easily seen to
lie inside $R^{\text{sym}}[\theta]P^{-n}$.
To conclude, we must apply the operator $E_\mathfrak{e}$ on both sides of (\ref{nexttolast}) after having multiplied 
by $\widetilde{\pi}$. On the left-hand side, 
we find
$$\sum_{i=1}^d\aa(t_i)(\mathcal{D}_{n}\omega)(t_i).$$
As for the right-hand side, the first term gives the function of $\TT_{t_1,\ldots,t_d}^\text{sym}$
$$E_\mathfrak{e}\left(\frac{\widetilde{\pi}(\mathcal{D}_{\theta,1}P)^nU_\aa}{P^{n+1}}\right).$$
The second and third terms, give $E_\mathfrak{e}(\Gamma_n+\Sigma_n)$ and yield an element $\Lambda_n$ of $$\sum_{i=1}^n\FF_q[t_1,\ldots,t_d]^{\text{sym}}E_\mathfrak{e}\left(\frac{\widetilde{\pi}\FF_q[\theta]}{P^i}\right),$$
and the Proposition follows.\CVD

\begin{Remark}{\em In this remark, not used in our paper, we point out an explicit way to compute the polynomials $U_\aa$
by specialization.
Let us consider indeterminates $X,X_1,\ldots,X_d,Y_1,\ldots,Y_d$ over a field $L$. Let us define the polynomials:
$$
P=P(X,X_1,\ldots,X_d)=\prod_{i=1}^d(X-X_d)
$$
and 
$$
U=U(X,X_1,\ldots,X_d,Y_1,\ldots,Y_d)=\sum_{i=0}^{d-1}(-1)^{d-1-i}\alpha_i^*(X_1,\ldots,X_d,Y_1,\ldots,Y_d)X^i,
$$
with 
$$
\alpha_{d-j}^*(X_1,\ldots,X_d,Y_1,\ldots,Y_d)=\sideset{{}^{}}{{}^{*}}\sum Y_{i_1}X_{i_2}\cdots X_{i_j},\quad j=1,\ldots,d
$$
where the sum $\sum^*$ is over the $j$-tuples $(i_1,i_2,\ldots,i_j)$ of pairwise distinct integers $1\leq i_1,\ldots,i_j\leq d$,
so that
\begin{eqnarray*}
\alpha_0^*&=&\sum_{j=1}^dY_j\prod_{i\neq j}X_i,\\
\alpha_{d-1}^*&=&\sum_{j=1}^dY_j.
\end{eqnarray*}

For example, if $d=3$, we have
\begin{eqnarray*}
\alpha_0^*&=&X_1X_2Y_3+X_1X_3Y_2+X_2X_3Y_1,\\
\alpha_1^*&=&X_1Y_2+X_1Y_3+X_2Y_1+X_2Y_3+X_3Y_1+X_3Y_2,\\
\alpha_2^*&=&Y_1+Y_2+Y_3,
\end{eqnarray*}
and $$U=\alpha_2^*X^2-\alpha_1^*X+\alpha_0^*.$$ If $\sigma$ is a permutation of $\{1,\ldots,d\}$,
then $$\alpha_j^*(X_{\sigma(1)},\ldots,X_{\sigma(d)},Y_{\sigma(1)},\ldots,Y_{\sigma(d)})=
\alpha_j(X_1,\ldots,X_d,Y_1,\ldots,Y_d)$$ for all $j$ (see \cite[Exercice 5, A IV.92]{Bou} for more information).
The following formula can be easily proved by induction
on $d\geq 1$:
$$
\sum_{i=1}^d\frac{Y_i}{X-X_i}=\frac{U}{P},$$
so that, for $\aa\in A$ and with the obvious choice of variables:
$$U_\aa(\theta,t_1,\ldots,t_d)=U(\theta,t_1,\ldots,t_d,\aa(t_1),\ldots,\aa(t_d)).$$}\end{Remark}

\subsection{Proof of Theorem \ref{mtre}}

Let $E^\infty$ be 
  the maximal abelian extension of $K$ tamely ramified at infinity which,
 by virtue of Hayes result in \cite{Ha1}, equals the field
 $$\FF_q^\text{alg}(\lambda_\mathfrak{a};a\in A).$$
 Furthermore, let $L^\infty$ be the field 
 $$\FF_q^\text{alg}((\mathcal{D}_n\omega)(\zeta);\zeta\in\FF_q^\text{alg},n\geq 0),$$ where 
 $\mathcal{D}_n$ denotes the $n$-th divided derivative with respect to the variable $t$ in $\TT_t$ (see Subsection \ref{chainrule}).
 Theorem \ref{mtre} states that $L^\infty=E^\infty$. Let $n$ be an integer and consider the set $\mathcal{E}_n$ whose 
 elements are the monic polynomials $\aa$ of $A$ such that for all $\pp$ a prime, $\pp^{n+1}$ does not divide $\aa$. Let us define the following subfields of $K^\text{alg}$:
 \begin{eqnarray*}
 E_n&=&\FF_q^\text{alg}(\phi_b(\lambda_\aa);\aa\in \mathcal{E}_n,b\in A)\\
 L_n&=&\FF_q^\text{alg}(\omega(\zeta),(\mathcal{D}_{1}\omega)(\zeta),\ldots,(\mathcal{D}_{n}\omega)(\zeta);\aa(\zeta)=0
 \text{ for some }\aa\in \mathcal{E}_{n+1}).
 \end{eqnarray*}
 We have $E_0\subset E_1\subset\cdots \subset E_n\subset\cdots$, 
 $L_0\subset L_1\subset\cdots \subset L_n\subset\cdots$, and $\cup_{i=0}^\infty E_i=E^\infty$, $\cup_{i=0}^\infty L_i=L^\infty$.
 We will prove (see Corollary \ref{acorollary})
 the identities $E_n=L_n$ by using induction on $n\geq 0$ as the next Proposition indicates; Theorem \ref{mtre} then follows 
 by a limit process.

The next result constitutes the main step to prove Theorem \ref{mtre}.
\begin{Proposition}\label{intermediateresult}
Let $\pp$ be a prime of degree $d$, let $\FF_\pp$ be the extension of $\FF_q$ generated by the roots
of $\pp$. For all $n\geq 0$, we have the following identity of fields.
$$\FF_\pp K_{\pp^{n+1}}=\FF_\pp((\mathcal{D}_k\omega)(\zeta),\zeta\in \FF_\pp,0\leq k\leq n).$$
\end{Proposition}
The proof of the above Proposition will proceed by induction on $n\geq0$. 
Let us denote by $L_{\pp,n}$ the field 
$\FF_\pp((\mathcal{D}_k\omega)(\zeta),\zeta\in \FF_\pp,0\leq k\leq n)$ and by $E_{\pp,n}$
the field $\FF_\pp K_{\pp^{n+1}}$. We have to show that $L_{\pp,n}=E_{\pp,n}$ for all $n\geq 0$. To ease the reading 
we will consider the case $n=0$ separately, before considering the general case, although this discrimination of the two cases is not
strictly necessary.

To prove the case $n=0$,
we will use the following Lemma.
\begin{Lemma}\label{lemmatrace}
Let $\pp$ be a prime of degree $d$, let $\zeta_1,\ldots,\zeta_d\in\FF_\pp$ be its roots. Let us also consider a polynomial $a\in A$. Then,
$$\sum_{j=1}^da(\zeta_j)\omega(\zeta_j)=\phi_{a\pp'}(\lambda_\pp).$$
\end{Lemma}
\noindent\emph{Proof.} Compare with \cite[Proposition 29]{gamma}. We recall that,
in the Tate algebra
$\TT_{t_1,\ldots,t_d}^{\text{sym}}$, we have, by the identity (\ref{secondfamily1}):
$$
\sum_{i=1}^da(t_i)\omega(t_i)=E_{\mathfrak{e}}\left(\frac{\widetilde{\pi}U_a}{P}\right),
$$
 with $P=\prod_{i=1}^d(\theta-t_i)$ and $U_a$ as in (\ref{primitiveformula}).
 If we replace $t_i=\zeta_i$ ($i=1,\ldots,d$) conjugate elements of $\FF_{q^d}$, we have $P(\theta,\zeta_1,\ldots,\zeta_d)=\pp$,
$U_a(\theta,\zeta_1,\ldots,\zeta_d)\equiv\pp' a\pmod{\pp}$ and the Lemma follows applying Lemma \ref{symmetric1} because:
$$\left.E_{\mathfrak{e}}\left(\frac{\widetilde{\pi}U_a}{P}\right)\right|_{t_i=\zeta_i}=\exp\left(\frac{\widetilde{\pi}\pp' a}{\pp}\right)=\phi_{a\pp'}(\lambda_\pp)$$ (remember that $U_a$ and $P$ are polynomials which are symmetric in $t_1,\ldots,t_d$).\CVD
\noindent\emph{Proof of Proposition \ref{intermediateresult} in the case $n=0$.} We begin by applying the first 
part of Proposition \ref{anidentity}, with $\aa=\pp$. In this case we get, for $\zeta\in\FF_\pp$ a root of $\pp$,
$\omega(\zeta)\in\FF_\pp(\omega_{\pp,0}^*(0),\ldots,\omega_{\pp,d-1}^*(0))=E_{\pp,0}$ (by (\ref{idcomposition})), from which the inclusion 
of fields $L_{\pp,0}\subset E_{\pp,0}$ follows. For the reverse inclusion we notice that $E_{\pp,0}=\FF_\pp(\lambda_\pp)=\FF_\pp(\phi_{\pp'}(\lambda_\pp))$
because $\pp$ and $\pp'$ are relatively prime. By Lemma \ref{lemmatrace} with $a=1$, we see that
$\phi_{\pp'}(\lambda_\pp)\in L_{\pp,0}$. Therefore, $L_{\pp,0}=E_{\pp,0}$.\CVD

To prove the case $n>0$ of the Proposition \ref{intermediateresult},
we will use two Lemmas.
The next Lemma is about the evaluation of the higher derivatives of
 the functions $\omega_{\aa,j}$ at roots of unity.
 \begin{Lemma}\label{lemmaevalzeta}
 Let $\pp$ be a prime of degree $d$ and $\zeta$ a root of $\pp$. Let $n\geq 1$ and $0\leq j\leq d-1$ be integers.
 Then, there exists an element $\mu_n\in E_{\pp,n-1}$ such that
 $$(\mathcal{D}_n\omega_{\pp,j})(\zeta)=(\pp'(\zeta))^n\phi_{\theta^j}(\lambda_{\pp^{n+1}})+\mu_n.$$
 \end{Lemma}
 \noindent\emph{Proof.} By (\ref{idcomposition}) and the chain rule, $(\mathcal{D}_n\omega_{\pp,j})(\zeta)$ is equal to
 $(\pp'(\zeta))^n(\mathcal{D}_{x,n}\omega^*_{\pp,j})(0)$ plus a linear combination of values $(\mathcal{D}_{x,i}\omega^*_{\pp,j})(0)$ with $i=0,\ldots,n-1$ and coefficients in $\FF_\pp$.\CVD

\begin{Lemma}\label{lemmatrace2}
Let $\pp$ be a prime of degree $d$, let $\zeta_1,\ldots,\zeta_d\in\FF_\pp$ be its roots, and let us consider an integer $n\geq 1$. Furthermore, let us also consider a polynomial $a\in A$. Then, there exist an element $\nu_n\in E_{\pp,n-1}$ such that 
$$\sum_{j=1}^da(\zeta_j)(\mathcal{D}_n\omega)(\zeta_j)=\phi_{a(\pp')^{n+1}}(\lambda_{\pp^{n+1}})+\nu_n.$$
\end{Lemma}
\noindent\emph{Proof.} By (\ref{secondfamily2}) and Lemma \ref{symmetric1}, we have
$$\sum_{i=1}^da(\zeta_i)(\mathcal{D}_n\omega)(\zeta_i)=\exp\left(\frac{(\pp')^{n+1}a}{\pp^{n+1}}\right)+\Lambda_n(\zeta_1,\ldots,\zeta_d).$$ The Lemma follows by setting 
$\nu_n=\Lambda_n(\zeta_1,\ldots,\zeta_d)$, belonging to $E_{\pp,n-1}=L_{\pp,n-1}$. Indeed, $\exp((\pp')^{n+1}a/\pp^{n+1})=\phi_{(\pp')^{n+1}a}(\lambda_{\pp^{n+1}})$.\CVD

\noindent\emph{Proof of Proposition \ref{intermediateresult} in the case $n\geq0$.} We proceed by induction on $n$;
the case $n=0$ being already proved, let us assume that $n>0$ and that $E_{\pp,i}=L_{\pp,i}$ for $i=0,\ldots,n-1$.
We first show that $L_{\pp,n}$ is contained in $E_{\pp,n}$. To do so, it suffices to show that 
$(\mathcal{D}_n\omega)(\zeta)\in E_{\pp,n}$ for all $\zeta$ root of $\pp$.

By (\ref{anidentityn}) with $t=\zeta$, we have 
$$(\mathcal{D}_n\omega)(\zeta)=\sum_{j=0}^{d-1}\pp_j(\zeta)(\mathcal{D}_n\omega_{\pp,j})(\zeta)+\Omega_n(\zeta),$$ 
where it is obvious that $\Omega_n(\zeta)\in L_{\pp,n-1}=E_{\pp,n-1}$. By Lemma \ref{lemmaevalzeta}, the sum
over $j$ equals
$$\sum_{j=0}^{d-1}\pp_j(\zeta)\pp'(\zeta)^n\phi_{\theta^j}(\lambda_{\pp^{n+1}})+\mu_n,$$ with $\mu_n\in E_{\pp,n-1}$.
This shows the inclusion $L_{\pp,n}\subset E_{\pp,n}$.

Let us now show the opposite inclusion. Since $\pp$ is prime, the class of $(\pp')^{n+1}$ is a generator of the $A$-module 
$A/\pp^{n+1}$. Therefore, in view of Lemma \ref{lemmatrace2}, we will only need to show that
$\phi_{(\pp')^{n+1}}(\lambda_{\pp^{n+1}})\in L_{\pp,n}$. But $\nu_n$ is element of $E_{\pp,n-1}=L_{\pp,n-1}$ 
by hypothesis, so that Lemma \ref{lemmatrace2} implies that $E_{\pp,n}\subset L_{\pp,n}$. The proof of Proposition
\ref{intermediateresult} is complete.\CVD

\begin{Corollary}\label{acorollary}
We have, for all $n\geq 0$, $E_n=L_n$.
\end{Corollary}
\noindent\emph{Proof}. On one side, $E_n$ is the compositum of all the fields $E_{\pp,n}$ with $\pp$ varying 
in the set of primes of $A$. On the other side, $L_n$ is the compositum of all the fields $L_{\pp,n}$, 
with $\pp$ varying in the set of primes of $A$. The result follows from Proposition \ref{intermediateresult}.\CVD

Theorem \ref{mtre} follows at once by taking $n\rightarrow\infty$.\CVD

\begin{Remark}{\em Theorem \ref{mtre} does not seem to be directly related to Anderson's result in \cite{An1}, where 
he proves that
the compositum of all the subfields of $\QQ^\text{alg}$ that are at once quadratic over $\QQ^{\text{ab}}$ and Galois over $\QQ$ is generated by {\em algebraic $\Gamma$-monomials} satisfying the Koblitz-Ogus condition. However, an analogue 
of our result concerning the Akhiezer-Baker function $\Gamma(s-t)$ as in \cite{gamma} and its composition with $t=$meromorphic integral-periodic functions should hold and will, we hope, be object of further investigation.}\end{Remark}
\subsection{Proof of Theorem \ref{gaussthakur}}
By the second part of Proposition \ref{Thakur}, we have that for $\delta\in \Delta_\pp$ and $\chi\in\widehat{\Delta}_\pp$,  
\begin{equation}\label{tha1}\delta (g (\chi ))=\chi (\delta) g (\chi).\end{equation}
Let $\zeta$ be the root of $\pp$ such that $\vartheta_\pp(\sigma_\theta)=\zeta$.
We have, for $j=0,\ldots,d-1$:
\begin{eqnarray*}
g(\vartheta_\pp^{q^j})&=&-\sum_{\delta\in\Delta_\pp}\vartheta_\pp(\delta^{-1})^{q^j}\delta(\lambda_\pp)\\
&=&-\sum_{a\in(A/\pp A)^\times}\vartheta_\pp(\sigma_a)^{-q^j}\sigma_a(\lambda_\pp)\\
&=&-\sum_{a\in(A/\pp A)^\times}\vartheta_\pp(\sigma_a)^{-q^j}\phi_a(\lambda_\pp).
\end{eqnarray*}
By (\ref{tha1}) and Lemma \ref{lemmatrace} ($\zeta_1,\ldots,\zeta_d$ again denote the zeros of $\pp$):
\begin{eqnarray*}
\chi_\zeta(\pp')^{q^j}g(\vartheta_\pp^{q^j})&=&\sigma_{\pp'}(g(\vartheta_\pp^{q^j}))\\
&=&-\sum_{a\in(A/\pp A)^\times}\vartheta_\pp(\sigma_{a\pp'})^{-q^j}\phi_{a\pp'}(\lambda_\pp)\\
&=&-\sum_{a\in(A/\pp A)^\times}\vartheta_\pp(\sigma_{a\pp'})^{-q^j}\sum_{k=1}^da(\zeta_k)\omega(\zeta_k)\\
&=&-\sum_{i=0}^{d-1}\omega(\zeta^{q^i})\sum_{a\in(A/\pp A)^\times}\vartheta_\pp(\sigma_{a})^{q^i-q^j},
\end{eqnarray*}
where the last identity follows from the fact that $\{\zeta_1,\ldots,\zeta_d\}=\{\zeta,\zeta^q,\ldots,\zeta^{q^{d-1}}\}$. Now, the sum
$$\sum_{a\in(A/\pp A)^\times}\vartheta_\pp(\sigma_{a})^{q^i-q^j}$$ always vanishes except when $i=j$, case in which the sum equals $-1$. Therefore:
$$\omega(\zeta^{q^j})=\omega(\vartheta_\pp(\sigma_\theta)^{q^j})=\chi_\zeta(\pp')^{q^j}g(\vartheta_\pp^{q^j})$$ hence completing the proof of the Theorem.\CVD

We shall also mention the following result.
\begin{Corollary}\label{corollary13}
Let $\zeta$ be the root of $\pp$ such that $\vartheta_\pp(\sigma_\theta)=\zeta$. The following identity holds: 
$$\FF_\pp K_\pp=\FF_\pp(g(\vartheta_\pp)).$$
\end{Corollary}
\noindent\emph{Proof.} By Proposition \ref{intermediateresult}, we have $\FF_\pp K_\pp=\FF_\pp(\omega(\zeta))$. 
Theorem \ref{gaussthakur} now implies that
$\FF_\pp(\omega(\zeta))=\FF_\pp(g(\vartheta_\pp)).$\CVD

\section{Functional identities for $L$-series\label{sectionthree}}

In this section, we prove Theorem \ref{functionalrelations}.
We will need a few preliminary results that we shall study in this subsection.
Let $d,s$ be non-negative integers. We begin with the study of the vanishing of the sums
$$S_{d,s}=S_{d,s}(t_1,\ldots,t_s)=\sum_{a\in A^+(d)}\chi_{t_1}(a)\cdots\chi_{t_s}(a)\in\FF_q[t_1,\ldots,t_s],$$
which are symmetric polynomials in $t_1,\ldots,t_s$ of total degree $\leq ds$, with the standard conventions on empty products.
We recall that, for $n\geq 0$, $$\sum_{a\in\FF_q}a^n$$ equals $-1$ if $n\equiv0\pmod{q-1}$ and $n\geq 1$, and equals 
$0$ otherwise.
We owe the next Lemma to D. Simon \cite{AngSim}. We give the proof here for the sake of completeness.
\begin{Lemma}[Simon's Lemma]\label{simon}
We have 
$S_{d,s}\neq0$ if and only if $d(q-1)\leq s$.
\end{Lemma}
\noindent\emph{Proof.} Since
$$S_{d,s}=\sum_{a_0 \in \FF_q}\cdots\sum_{a_{d-1} \in \FF_q}\prod_{i=1}^s (a_0+a_1t_i+\cdots + a_{d-1} t_i^{d-1}+t_i^d),$$
the coefficient $c_{v_1,\ldots,v_s}$ of $t_1^{v_1}\cdots t_s^{v_s}$ with $v_i\leq d$ ($i=1,\ldots,s$) is given by the sum:
$$\sum_{a_0 \in \FF_q}\cdots\sum_{a_{d-1} \in \FF_q}a_{v_1}\cdots a_{v_s},$$ if we set $a_{d}=1$.
The last sum can be rewritten as:
\begin{equation}\label{c1cs}
c_{v_1,\ldots,v_s}=\left(\sum_{a_0\in\FF_q}a_0^{\mu_0}\right)\cdots\left(\sum_{a_{d-1}\in\FF_q}a_{d-1}^{\mu_{d-1}}\right),
\end{equation}
where $\mu_i$ is the cardinality of the set of the indices $j$ such that $v_j=i$, from which one notices that $$\sum_{i=0}^{d-1}\mu_i\leq s$$
(notice also that $s-\sum_i\mu_i$ is the cardinality of the set of indices $j$ such that $v_j=d$).
For any choice of $\mu_0,\ldots,\mu_{d-1}$ such that $\sum_i\mu_i\leq s$, there exists $(v_1,\ldots,v_s)$ such that
(\ref{c1cs}) holds. 

If $s<d(q-1)$, for all $(v_1,\ldots,v_s)$ as above, there exists $i$ such that, in (\ref{c1cs}), $\mu_i<q-1$ so that $S_{d,s}=0$.
On the other hand, if $s\geq d(q-1)$, it is certainly possible to find $(v_1,\ldots,v_s)$ such that, in (\ref{c1cs}),
$\mu_0=\cdots=\mu_{d-1}=q-1$ so that the sum does not vanish in this case.\CVD

As an immediate corollary of Lemma \ref{simon}, we see that the series
$$F_s=F_s(t_1,\ldots,t_s)=\sum_{d\geq0}S_{d,s}=\sum_{d\geq 0}\sum_{a\in A^+(d)}\chi_{t_1}(a)\cdots\chi_{t_s}(a)$$
defines a symmetric polynomial of $\FF_q[t_1,\ldots,t_s]$ of total degree at most $\frac{s^2}{q-1}$.
In the next Lemma, we provide a sufficient condition for the vanishing of the polynomial $F_s$.
\begin{Lemma}\label{vanishingFs} If $s\geq 1$, then, $F_s=0$ if and only if 
$s\equiv0\pmod{q-1}$.
\end{Lemma}
\noindent\emph{Proof.} Let us assume first that 
$s\equiv0\pmod{q-1}$. The hypothesis on $s$ implies that
$$\sum_{a\in A,\deg_\theta(a)=d}\chi_{t_1}(a)\cdots\chi_{t_s}(a)=-S_{d,s}.$$
We denote by $A(\leq d)$ the set of polynomials of $A$ of degree $\leq d$ and we write
$$G_{d,s}=\sum_{a\in A(\leq d)}\chi_{t_1}(a)\cdots\chi_{t_s}(a).$$
We then have:
$$G_{\frac{s}{q-1},s}=-F_s.$$
Let us choose now distinct primes $\pp_1,\ldots,\pp_s$ of respective degrees $d_1,\ldots,d_s\geq s/(q-1)$ and $\mathfrak{f}=\pp_1\cdots\pp_s$. 
For all $i=1,\ldots,s$, we choose a root $\zeta_i\in\FF_q^{\text{alg}}$ of $\pp_i$. 
Let us then consider the Dirichlet character of the first kind $\chi=\chi_{\zeta_1}\cdots\chi_{\zeta_s}$.
We have:
\begin{eqnarray*}
F_s(\zeta_1,\ldots,\zeta_s)&=&-G_{\frac{s}{q-1},s}(\zeta_1,\ldots,\zeta_s)\\
&=&-\sum_{a\in A(\leq s/(q-1))}\chi(a)\\
&=&-\sum_{a\in A(\leq d_1+\cdots+d_s)}\chi(a)\\
&=&-\sum_{a\in (A/\mathfrak{f}A)^\times}\chi(a)\\
&=&0,
\end{eqnarray*}
by \cite[Proposition 15.3]{ROS}. Since the set  of $s$-tuples $(\zeta_1,\ldots,\zeta_s)\in(\FF_{q}^{\text{alg}})^s$
with $\zeta_1,\ldots,\zeta_s$ as above is Zariski-dense in $\mathbb{A}^s(\CC_\infty)$, this implies the vanishing of $F_s$.
On the other hand, if $s\not\equiv0\pmod{q-1}$, then $F_s(\theta,\ldots,\theta)=\zeta(-s)$ the $s$-th ``odd negative" Goss' zeta value which is non-zero,
see \cite{Go0}.

\CVD

\subsection{Analyticity\label{analyticity}}

The functions $L(\chi_{t_1}\cdots\chi_{t_s},\alpha)$ are in fact rigid analytic entire functions of $s$ variables.
This property, mentioned in \cite{Pe}, can be deduced from the more general 
Proposition \ref{analytic} that we give here for convenience of the reader.

Let $a$ be a monic polynomial  of $A$. we set:
$$\langle a\rangle=\frac{a}{\theta^{\deg_\theta(a)}}\in 1+\theta^{-1}\mathbb{F}_q[\theta^{-1}].$$
Let $y\in\mathbb{Z}_p$, where $p$ is the prime dividing $q$. Since $\langle a\rangle$ is a $1$-unit of $K_\infty$, 
we can consider its exponentiation by $y$:
$$\langle a\rangle^{y}=\sum_{j\geq0}\binom{y}{j}(\langle a\rangle-1)^j\in \FF_q[[\theta^{-1}]].$$
Here, the binomial $\binom{y}{j}$ is defined, for $j$ a non-negative integer, by extending Lucas formula:
writing the $p$-adic expansion $\sum_{i\geq 0}y_ip^i$ of $y$ ($y_i\in\{0,\ldots,p-1\}$) and the $p$-adic expansion $\sum_{i=0}^rj_ip^i$ of $j$ ($j_i\in\{0,\ldots,p-1\}$),
we are explicitly setting:
$$\binom{y}{j}=\prod_{i=0}^r\binom{y_i}{j_i}.$$

We also recall, from \cite[Chapter 8]{Go2}, the topological group $\mathbb{S}_\infty=\CC_\infty^\times\times\ZZ_p$.
For $(x,y)\in\mathbb{S}_\infty$ and $d,s$ non-negative integers, we define the sum
$$S_{d,s}(x,y)=S_{d,s}(x,y)(t_1,\ldots,t_s)=x^{-d}\sum_{a\in A^+(d)}\chi_{t_1}(a)\cdots\chi_{t_s}(a)\langle a \rangle^y\in x^{-d}K_\infty[t_1,\ldots,t_s],$$
which is, for all $x,y$, a symmetric polynomial of total degree $\leq ds$.

Let us further define, more generally, for variables $t_1,\ldots,t_s\in\CC_\infty$ and $(x,y)\in\mathbb{S}_\infty$, the series:
$$L(\chi_{t_1}\ldots\chi_{t_s};x,y)=\sum_{d\geq0}S_{d,s}(x,y)(t_1,\ldots,t_s).$$ For fixed choices of $(x,y)\in\mathbb{S}_\infty$,
it is easy to show that
$$L(\chi_{t_1}\ldots\chi_{t_s};x,y)\in\CC_\infty[[t_1,\ldots,t_s]],$$ and with a little additional work, one also verifies
that this series defines an element of $\TT_{t_1,\ldots,t_s}$. Of course, if $(x,y)=(\theta^\alpha,-\alpha)$ with $\alpha>0$ integer, we find
$$L(\chi_{t_1}\ldots\chi_{t_s};\theta^\alpha,-\alpha)=L(\chi_{t_1}\ldots\chi_{t_s},\alpha).$$

The next Proposition holds and improves results of Goss; see \cite[Theorems 1, 2]{Go3}).

\begin{Proposition}\label{analytic}
The series $L(\chi_{t_1},\ldots,\chi_{t_s};x,y)$ converges for all $(t_1,\ldots,t_s)$
and for all $(x,y)\in\mathbb{S}_\infty$, to a 
continuous-analytic function on $\CC_\infty^s\times\mathbb{S}_\infty$ in the sense of Goss.
\end{Proposition}

The proof of this result is a simple consequence of the Lemma below.
The norm $\|\cdot\|$ used in the Lemma is that of $\TT_{t_1,\ldots,t_s}$.

\begin{Lemma}
\label{lemma_bis}
Let $(x,y)$ be in $\mathbb{S}_\infty$ and let us consider an integer $d> (s+1)/(q-1)$, with $s>0$. Then:
$$\|S_{d,s}(x,y)\|\leq|x|^{-d}q^{-q^{\lfloor d-\frac{s+1}{q-1}\rfloor}}.$$
\end{Lemma}
\noindent\emph{Proof.} 
Let us write the $p$-adic expansion $y=\sum_{n\geq 0}c_n p^n,$ with $c_n \in \{ 0,\ldots , p-1\}$ for all $n$. 
Collecting blocks of $e$ consecutive terms (where $q=p^e$), this yields a ``$q$-adic" expansion, from which we can extract partial sums:
$$
y_n=\sum_{k=0}^{en-1}c_kp^k=\sum_{i=0}^{n-1}u_iq^i\in\ZZ_{\geq 0},
$$
where $$u_i=\sum_{j=ei}^{e(i+1)-1}c_jp^{j-ei}\in\{0,\ldots,q-1\}.$$
In particular, for $n\geq 0,$ we observe that
$\ell_q(y_n)\leq n(q-1).$
Since
\begin{eqnarray*}
S_{d,s}(x,y_n)&=&\frac{1}{x^d\theta^{dy_n}}\sum_{a\in A^+(d)} \chi_{t_1}(a)\cdots\chi_{t_s}(a)a^{y_n}\\
&=&\frac{1}{x^d\theta^{dy_n}}S_{d,r}(t_1,\ldots,t_s,\underbrace{\theta,\ldots,\theta}_{u_0\text{ times }},\underbrace{\theta^q,\ldots,\theta^q}_{u_1\text{ times }},\ldots,\underbrace{\theta^{q^{n-1}},\ldots,\theta^{q^{n-1}}}_{u_{n-1}\text{ times }})\end{eqnarray*} with $r=s+\ell_q(y_n)$,
if $d(q-1)>s+\ell_q(y_n),$ we have by Simon's Lemma \ref{simon}:
$$S_{d,s}(x,y_n)=0.$$ This condition is ensured if $d(q-1)>s+n(q-1)$.

Now, we claim that
$$\|S_{d,s}(x,y)-S_{d,s}(x,y_n)\|\leq |x|^{-d}q^{-q^{n}}.$$
Indeed,
$$S_{d,s}(x,y)-S_{d,s}(x,y_n)=x^{-d}\sum_{a\in A^+(d)}\chi_{t_1}(a)\cdots\chi_{t_s}(a)
\sum_{j\geq0}\left(\binom{y}{j}-\binom{y_n}{j}\right)(\langle a\rangle-1)^j,$$
and $\binom{y}{j}=\binom{y_n}{j}$ for $j=0,\ldots,q^{n}-1$ by Lucas' formula and the definition of the binomial, so that
$$\left|\sum_{j\geq0}\left(\binom{y}{j}-\binom{y_n}{j}\right)(\langle a\rangle-1)^j\right|\leq q^{-q^{n}}.$$

The Lemma follows by choosing $n=\lfloor d-1-\frac{s+1}{q-1}\rfloor$.\CVD

In particular, we have the following Corollary to Proposition \ref{analytic} which generalizes 
\cite[Theorem 1]{Go3}, the deduction of which, easy, is left to the reader.
\begin{Corollary}\label{corollaryanalytic}
For any choice of an integer $\alpha>0$ and non-negative integers $M_1,\ldots,M_s$, the function
$$L(\chi_{t_1}^{M_1}\cdots\chi_{t_s}^{M_s},\alpha)=\sum_{d\geq 0}\sum_{a\in A^+(d)}\chi_{t_1}(a)^{M_1}\cdots\chi_{t_s}(a)^{M_s}a^{-\alpha}$$
defines a rigid analytic entire function $\CC_\infty^s\rightarrow\CC_\infty$.
\end{Corollary}

\subsection{Computation of polynomials with coefficients in $K_\infty$}

\begin{Lemma}\label{lemmaremark}
For all $d\geq 0$, we have:
$$S_d(-\alpha)=\sum_{a\in A^+(d)}a^{-\alpha}\not = 0.$$
\end{Lemma}
\noindent\emph{Proof.} This follows from \cite[proof of Lemma 8.24.13]{Go2}.\CVD

We introduce, for $d,s,\alpha$ nonnegative integers, the sum:
$$S_{d,s}(-\alpha)=\sum_{a\in A^+(d)}\chi_{t_1}(s)\cdots\chi_{t_s}(a)a^{-\alpha}\in K[t_1,\ldots,t_s],$$ representing a symmetric polynomial
of $K[t_1,\ldots,t_s]$ of exact total degree $ds$ by Lemma \ref{lemmaremark}. We have, with the notations of Section \ref{analyticity}:
$$S_{d,s}(-\alpha)=S_{d,s}(\theta^\alpha,-\alpha).$$
From the above results, we deduce the following Proposition.
\begin{Proposition}\label{firstintermprop}
Let $l\geq 0$ be an integer such that $q^{l}-\alpha\geq 0$ and $2\leq  \ell_q(q^{l}-\alpha)+s\leq d(q-1)$. Then:
$$S_{d,s}(-\alpha) \equiv 0\pmod{\prod_{j=1}^s(t_j-\theta^{q^{l}})}.$$
Furthermore, assume that
$s\equiv \alpha \pmod{q-1}.$ With $l$ as above, let $k$ be an integer such that $k(q-1)\geq \ell_q(q^l-\alpha)+s.$ Then:
$$\sum_{d=0}^{k} S_{d,s}(-\alpha)\equiv 0\pmod{\prod_{j=1}^s(t_j-\theta^{q^{l}})}.$$
\end{Proposition}
\noindent\emph{Proof.} Let us write $m=\ell_q(q^{l}-\alpha)$. We have 
$s-1+m<d(q-1)$ so that, by Simon's Lemma \ref{simon}, $S_{d,s-1+m}=0$.
Now, let us write the $q$-ary expansion $q^l-\alpha=n_0+n_1q+\cdots+n_rq^r$ with $n_i\in\{0,\ldots,q-1\}$ and let us observe that, since $q^l-\alpha\geq0$, 
\begin{eqnarray*}
S_{d,s}(-\alpha)(t_1,\ldots,t_{s-1},\theta^{q^l})&=&\sum_{a\in A^+(d)}\chi_{t_1}(a)\cdots\chi_{t_{s-1}}(a)a^{q^l-\alpha}\\
&=&\sum_{a\in A^+(d)}\chi_{t_1}(a)\cdots\chi_{t_{s-1}}(a)\chi_{\theta}(a)^{n_0}\chi_{\theta^q}(a)^{n_1}\cdots\chi_{\theta^{q^r}}(a)^{n_r}\\
&=&S_{d,s-1+m}(t_1,\ldots,t_{s-1},\underbrace{\theta,\ldots,\theta}_{n_0\text{ times }},\underbrace{\theta^q,\ldots,\theta^q}_{n_1\text{ times }},\ldots,\underbrace{\theta^{q^r},\ldots,\theta^{q^r}}_{n_r\text{ times }})\\
&=&0.
\end{eqnarray*}
Therefore $t_s-\theta^{q^{l}}$ divides $S_{d,s}(-\alpha).$ The first part of the Proposition follows from the fact that this polynomial is symmetric.
For the second part, we notice by the first part, that the condition on $k$ is sufficient for the sum $S_{d,s}(-\alpha)(t_1,\ldots,t_{s})$ to be 
congruent modulo $(t_s-\theta^{q^l})$ for all
$d\geq k+1$. On the other hand, by Lemma \ref{vanishingFs} and the above computation, we have
$$\sum_{d\geq0}S_{d,s}(-\alpha)(t_1,\ldots,t_{s-1},\theta^{q^l})=F_s(t_1,\ldots,t_{s-1},\theta,\ldots,\theta^{q^r})=0.$$
But then, thanks to the condition on $k$,
$$\sum_{d=0}^kS_{d,s}(-\alpha)\equiv-\sum_{d>k}S_{d,s}(-\alpha)\equiv0\pmod{(t_s-\theta^{q^l})}$$ and the Proposition follows again because 
the sum we are inspecting is a symmetric polynomial.\CVD 

We further have the result below.

\begin{Proposition}
\label{proposition1}
Let $s,\alpha\geq 1, $ $s\equiv \alpha\pmod{q-1}.$  Let $\delta$ be the smallest positive integer such that 
$q^{\delta}\geq\alpha$ and $s+  \ell_q(q^{\delta}-\alpha)\geq 2.$ Then, the function of Theorem \ref{functionalrelations}
$$V_{\alpha,s}(t_1,\ldots,t_s)=L(\chi_{t_1}\cdots\chi_{t_s},\alpha)\omega(t_1)\cdots\omega(t_s)
\widetilde{\pi}^{-\alpha}\left(\prod_{i=1}^{s}\prod_{j=0}^{\delta-1}\left(1-\frac{t_i}{\theta^{q^j}}\right)\right)$$
is in fact a symmetric polynomial of $K_{\infty}[t_1,\ldots ,t_s]$. Moreover, its 
total degree $\delta(\alpha,s)$ is not bigger than $s\left(\frac{s+\ell_q(q^{\delta}-\alpha)}{q-1}\right)-s.$
\end{Proposition}
\noindent\emph{Proof.} Let $\delta$ be the smallest positive integer such that $q^{\delta}-\alpha\geq 0$ and $s+  \ell_q(q^{\delta}-\alpha)\geq 2.$ We fix  an integer $k$ such that  \begin{equation}\label{withdelta}k(q-1)\geq s+\ell_q(q^{\delta}-\alpha).\end{equation}
We also set:
$$N(k)=\delta+k-\frac{s+\ell_q(q^\delta-\alpha)}{q-1}.$$
Obviously, $N(k)\geq \delta$.
Let $l$ be an integer such that 
$$
\delta\leq l\leq N(k).
$$
We claim that we also have
$$k(q-1)\geq s+\ell_q(q^{l}-\alpha).$$
Indeed,
let us write the $q$-ary expansion $\alpha=\alpha_0+\alpha_1 q+\cdots+\alpha_mq^m$ with $\alpha_m\neq0$.
Then, $\delta=m$ if $\alpha=q^m$ and $s\geq 2$ and $\delta=m+1$ otherwise. If $l$ is now an integer $l\geq \delta$,
we have
\begin{eqnarray*}
q^l-\alpha&=&q^l-q^{\delta}+q^{\delta}-\alpha\\
&=&q^{\delta}(q-1)\left(\sum_{i=0}^{l-\delta-1}q^i\right)+q^\delta-\alpha,
\end{eqnarray*}
where the sum over $i$ is zero if $l=\delta$, and
$$\ell_q(q^l-\alpha)=(q-1)(l-\delta)+\ell_q(q^\delta-\alpha)$$ because there is no carry over in the above sum. Now, 
the claim follows from (\ref{withdelta}).

By Proposition \ref{firstintermprop} we have, with $k$ as above, that the following expression
$$W_{k,s,\alpha}:=\left(\prod_{i=1}^{s}\prod_{j=\delta}^{N(k)}\left(1-\frac{t_i}{\theta^{q^j}}\right)^{-1}\right)\sum_{d=0}^{k} S_{d,s}(-\alpha)$$
is in fact a symmetric polynomial in $K[t_1,\ldots , t_s]$.
By Lemma \ref{lemmaremark}, $S_{d,s}(-\alpha)\in K[t_1, \ldots ,t_s]$ is a symmetric polynomial of total degree $ds$; indeed, 
the coefficient of $t_1^d\cdots t_s^d$ is exactly $S_d(-\alpha).$ Hence, the total degree of $\sum_{d=0}^{k} S_{d,s}(-\alpha)$
is exactly $ks$. The total degree of the product $$\prod_{i=1}^{s}\prod_{j=\delta}^{N(k)}\left(1-\frac{t_i}{\theta^{q^j}}\right)$$
is equal to $s(1+N(k)-\delta)$ so that, by the definition of $N(k)$:
\begin{eqnarray*}
\deg(W_{k,s,\alpha})&=&sk-s-sN(k)+s\delta\\
&=&sk-sk-s\delta+s\delta-s +s\left(\frac{s+\ell_q(q^\delta-\alpha)}{q-1}\right)\\
&=&s\left(\frac{s+\ell_q(q^\delta-\alpha)}{q-1}\right)-s,
\end{eqnarray*}
independent on $k$. We now let $k$ tend to infinity. The Proposition follows directly from the definition (\ref{definitionomega})
of $\omega$ as an infinite product, the fact that, in (\ref{definitionomega}), $\widetilde{\pi}\theta^{-1/(q-1)}\in K_\infty$, and  the definition of 
$L(\chi_{t_1}\cdots\chi_{t_s},\alpha)$.\CVD

\subsection{An intermediate result on special values of Goss $L$-functions}

Let $\chi$ be a Dirichlet character of the first kind, that is, a character $$\chi :(A/\aa A)^\times \rightarrow (\FF_q^{\text{alg}})^\times,$$ where $\aa$ is 
a non-constant squarefree monic element of $A$ which we identify, by abuse of notation, to a character of $\widehat{\Delta}_\aa$ still denoted by $\chi$, of conductor $\mathfrak{f}=\mathfrak{f}_\chi$, and degree $d=\deg_\theta f$. 

Let $s(\chi)$ be the {\em type} of $\chi$, that is, the unique integer $s(\chi)\in \{ 0,\ldots , q-2\}$ such that:
$$\chi(\zeta)=\zeta^{s(\chi)}\quad \text{ for all }\zeta\in \FF_q^\times.$$
We now consider the {\em generalized $\alpha$-th Bernoulli number} $B_{\alpha,\chi^{-1}}\in \mathbb F_q(\chi )(\theta)$ associated to $\chi^{-1},$ 
\cite[Section 2]{ANG&TAE}, and the special value of Goss' {\em abelian $L$-function} \cite[Section 8]{Go2}: $$L(\alpha,\chi)=\sum_{a\in A^+}\chi(a)a^{-\alpha},\quad \alpha\geq 1.$$
The following result is inspired by the proof of \cite[Proposition 8.2]{ANG&TAE}:
\begin{Proposition}
\label{Lvalues}  Let $\alpha\geq 1, $ $\alpha\equiv s(\chi) \pmod{q-1}.$ Then:
$$\frac{L(\alpha,\chi)g(\chi)}{\widetilde{\pi}^\alpha}= (-1)^{d}\frac{B_{\alpha,\chi^{-1}}}{\mathfrak{f}^{\alpha-1}}\in\FF_q(\chi)(\theta).$$
\end{Proposition}
\noindent\emph{Proof.} 
The proposition is known to be true for the trivial character (see \cite[Section 9.2]{Go2}); in this case, we notice that:
$$B_{\alpha,\chi_0^{-1}}=\frac{BC_\alpha}{\Pi (\alpha)},\quad \alpha\geq 1,\quad \alpha\equiv 0\pmod{q-1},$$
where we recall that $BC_\alpha$ is the $\alpha$-th Bernoulli-Carlitz number and $\Pi (\alpha)$  is the Carlitz factorial of $\alpha$
(see \cite[Definition 9.2.1]{Go2}). 
We now assume that $\chi \not = \chi_0$. Since:
 $$\exp(z)=z\prod_{a\in A\setminus \{0\}}\left(1-\frac{z}{\widetilde{\pi}a}\right),$$
 We have:
 $$\frac{1}{\exp(z)}=\sum_{a\in A}\frac{1}{z-\widetilde{\pi}a}.$$
 Let $b\in A$ be relatively prime with $\mathfrak{f}$ and let $\sigma_b \in \mathbf{Gal}(K_\mathfrak{f}/K)$
 be the element such that $\sigma_b(\lambda_\mathfrak{f})= \phi_b(\lambda_\mathfrak{f}).$ We have:
 $$\frac{1}{\exp(z)-\sigma_b(\lambda_\mathfrak{f})}  =-\sum_{n\geq 0} \frac{\mathfrak{f}^{n+1}}{\widetilde{\pi}^{n+1}}\left(\sum_{a \in A}\frac{1}{(b+a\mathfrak{f})^{n+1}}\right)z^n.$$
Therefore,   we obtain:
 $$\sum_{b\in (A/\mathfrak{f}A)^\times} \frac{\chi (b)}{\exp(z)-\sigma_b(\lambda_\mathfrak{f})} = -\sum_{n\geq 0}\frac{\mathfrak{f}^{n+1}}{\widetilde{\pi}^{n+1}}\left(\sum_{a\in A\setminus\{0\}}\frac{\chi (a)}{a^{n+1}}\right)z^n.$$
 If $n+1\not \equiv s(\chi )\pmod{q-1},$ we get:
 $$\sum_{a\in A\setminus\{0\}}\frac{\chi (a)}{a^{n+1}}=0,$$
 and if $n+1\equiv s(\chi )\pmod{q-1},$ we have:
 $$\sum_{a\in A\setminus\{0\}}\frac{\chi (a)}{a^{n+1}}=-L(n+1, \chi).$$
 Thus:
 \begin{equation}\label{identitytobeusedagain}
 \sum_{b\in (A/\mathfrak{f}A)^\times} \frac{\chi (b)}{\exp(z)-\sigma_b(\lambda_\mathfrak{f})} = \sum_{i\geq 1,\, i\equiv s(\chi)\pmod{(q-1)}}\frac{\mathfrak{f}^{i}L(i,\chi)}{\widetilde{\pi}^i}z^{i-1}.\end{equation}
 But note that by  the second part of Lemma \ref{normalbasis}:
 $$\sum_{b\in (A/\mathfrak{f}A)^\times} \frac{\chi (b)}{\exp(z)-\sigma_b(\lambda_\mathfrak{f})} \in g(\chi^{-1})\mathbb \mathfrak{F}_q(\chi)(\theta)[[z]].$$
Since by Proposition \ref{Thakur}, $$g(\chi)g(\chi^{-1} )= (-1)^d\mathfrak{f},$$ where $d=\deg_\theta f_\chi$, we get the result by comparison of the coefficients 
of the series expansion of both sides of (\ref{identitytobeusedagain}).\CVD

\begin{Remark}\label{remarktobeusedlater}{\em In the above proof of Proposition \ref{Lvalues}, if we set $\alpha=1$ we have, by comparison of the constant 
terms in the series expansions in powers of $z$ in (\ref{identitytobeusedagain}):
$$\widetilde{\pi}^{-1}\mathfrak{f}L(1,\chi)=-\sum_{b\in(A/\mathfrak{f}A)^\times}\frac{\chi(b)}{\sigma_b(\lambda_\mathfrak{f})}\in g(\chi^{-1})\FF_q(\chi)(\theta).$$
Assuming that $\mathfrak{f}$ is not a prime, by \cite[Proposition 12.6]{ROS}, $\lambda_\mathfrak{f}$ is a unit in the integral closure $ A_{\mathfrak{f}}$
of $A$ in $K_{\mathfrak{f}}$. Therefore, $$\sum_{b\in(A/\mathfrak{f}A)^\times}\frac{\chi(b)}{\sigma_b(\lambda_\mathfrak{f})}\in g(\chi^{-1})\FF_q(\chi)[\theta]
$$ and we deduce that
$$\widetilde{\pi}^{-1}L(1,\chi)g(\chi)\in\FF_q(\chi)[\theta].$$ This remark will be crucial in the proof of Corollary \ref{corollaryL}.}\end{Remark}

\subsection{Proof of Theorem \ref{functionalrelations}}

The next Lemma provides a rationality criterion for a polynomial a priori with coefficients in $K_\infty$,
again based on evaluation at roots of unity.
\begin{Lemma}
\label{lemma3} 
Let $F(t_1,\ldots, t_s) \in K_{\infty}[t_1,\ldots ,t_s]$ such that for all $\zeta_1,\ldots ,\zeta_s \in\FF_q^{\text{alg}},$ pairwise not conjugate over $\FF_q$, $$F(\zeta_1,\ldots ,\zeta_s)\in K(\zeta_1,\ldots ,\zeta_s).$$ Then $F(t_1,\ldots ,t_s)\in K[t_1,\ldots ,t_s].$ 
\end{Lemma}
\noindent\emph{Proof.} 
We begin by pointing out that if elements $a_1,\ldots,a_s\in K_\infty$ are $K\otimes_{\FF_q}\FF_q^{\text{alg}}$-linearly dependent, then they also 
are $K$-linearly dependent. The proof proceeds by induction on $s\geq 1$. For $s=1$, this is obvious. Now, let 
\begin{equation}\label{test1}
\sum_{i=1}^s\lambda_sa_s=0
\end{equation} be a non-trivial relation of linear dependence with the $\lambda_i\in K\otimes\FF_q^{\text{alg}}\setminus\{0\}$.
We may assume that $\lambda_s=1$ and that there exists $i\in\{1,\ldots,s-1\}$ such that 
$\lambda_i\not\in K$. Then, there exists $$\sigma\in\mathbf{Gal}(K_\infty\otimes\FF_q^{\text{alg}}/K_\infty)=
\mathbf{Gal}(K\otimes\FF_q^{\text{alg}}/K)=\mathbf{Gal}(\FF_q^{\text{alg}}/\FF_q)$$ such that
$\sigma(\lambda_i)\neq\lambda_i$. Applying $\sigma$ on both left- and right-hand sides of (\ref{test1}) and subtracting,
yields a non-trivial relation involving at most $s-1$ elements of $K_\infty$ on which we can apply the induction hypothesis.

We can now complete the proof of the Lemma. Let 
$F$ be a polynomial in $K_\infty[t_1,\ldots,t_s]$ not in $K[t_1,\ldots,t_s]$. It is easy to show that there exist
$a_1,\ldots,a_m\in K_\infty$, linearly independent over $K$, such that 
$$F=a_1P_1+\cdots+a_mP_m,$$ where $P_1,\ldots,P_m$ are non-zero polynomials of $K[t_1,\ldots,t_s]$.
Let us suppose by contradiction that there exists $F\in K_\infty[t_1,\ldots,t_s]\setminus K[t_1,\ldots,t_s]$ satisfying the hypotheses of the Lemma. Since the 
set of $s$-tuples $(\zeta_1,\ldots,\zeta_s)$ as in the statement of the Lemma is Zariski-dense in $\mathbb{A}^s(\CC_\infty)$, there exist a choice of such roots of unit $\zeta_1,\ldots,\zeta_s$ and $i\in\{1,\ldots,m\}$ such that $P_i(\zeta_1,\ldots,\zeta_m)\neq0.$ This means that $a_1,\ldots,a_m$ are $K\otimes\FF_q^{\text{alg}}$-linearly dependent, thus $K$-linearly dependent by the previous observations;
a contradiction.\CVD

\noindent\emph{Proof of Theorem \ref{functionalrelations}.} In view of Lemma \ref{lemma3}, we want to show that the polynomial

$$V_{\alpha,s}=\widetilde{\pi}^{-\alpha}L(\chi_{t_1}\cdots\chi_{t_s},\alpha)\omega(t_1)\cdots\omega(t_s)
\left(\prod_{i=1}^{s}\prod_{j=0}^{\delta-1}\left(1-\frac{t_i}{\theta^{q^j}}\right)\right)
\in K_{\infty}[t_1,\ldots ,t_s]$$ of Proposition \ref{proposition1} takes values in $K(\zeta_1,\ldots ,\zeta_s)$ for 
all $\zeta_1,\ldots ,\zeta_s \in\FF_q^{\text{alg}}$ pairwise non conjugate over $\FF_q$. Let $(\zeta_1,\ldots,\zeta_s)$ be one of such 
$s$-tuples of roots of unity and, for $i=1,\ldots,s$, let $\pp_i\in A$ be the minimal polynomial of $\zeta_i$, so that $\pp_1,\ldots,\pp_s$ are
pairwise relatively prime. We choose the characters $\vartheta_{\pp_i}$ so that
$\vartheta_{\pp_i}(\sigma_\theta)=\zeta_i$ for all $i$. We construct the Dirichlet character of the first kind $\chi$ defined, for $a\in A$, by $$\chi(a)=\chi_{\zeta_1}(a)\cdots\chi_{\zeta_s}(a).$$
By Proposition \ref{Lvalues}, we have
$$\frac{L(\alpha,\chi)g(\chi)}{\widetilde{\pi}^\alpha}=(-1)^{d_\chi}\frac{B_{\alpha,\chi^{-1}}}{\mathfrak{f}_{\chi}^{\alpha-1}}\in \FF_q(\chi)(\theta).$$
Since
$$L(\alpha,\chi)=L(\chi_{\zeta_1}\cdots\chi_{\zeta_s},\alpha),$$ we get:
\begin{eqnarray*}
V_{\alpha,s}(\zeta_1,\ldots,\zeta_s)&=&L(\alpha,\chi)\omega(\zeta_1)\cdots\omega(\zeta_s)\widetilde{\pi}^{-\alpha}\\
&=&\frac{L(\alpha,\chi)g(\chi)}{\widetilde{\pi}^{\alpha}}\frac{\omega(\zeta_1)\cdots\omega(\zeta_s)}{g(\chi)}\\
&=&(-1)^{d_\chi}\frac{B_{\alpha,\chi^{-1}}}{\mathfrak{f}_{\chi}^{\alpha-1}}\chi_{\zeta_1}(\pp_1')\cdots\chi_{\zeta_s}(\pp_s')\\
&\in& K(\zeta_1,\ldots,\zeta_s),
\end{eqnarray*}
where in the next to last step, we have used Theorem \ref{gaussthakur}. The proof of Theorem \ref{functionalrelations} now follows from Lemma \ref{lemma3}.\CVD

\section{Congruences for Bernoulli-Carlitz numbers\label{congruences}}

In this Section, we shall prove Theorem \ref{bernoulli}. This is possible because 
in Theorem \ref{functionalrelations}, more can be said when $\alpha=1$.
In this case, one sees that 
the integer $\delta$ of Theorem \ref{functionalrelations} is equal to zero and $s\geq q$, so that, with the notations of that result,
$$V_{1,s}=\widetilde{\pi}^{-1}L(\chi_{t_1}\cdots\chi_{t_s},1)\omega(t_1)\cdots\omega(t_s).$$
In the next Subsection we will show that this is a polynomial of $A[t_1,\ldots,t_s]$.

\subsection{Functional identities with $\alpha=1$}

We begin with the following Corollary of Theorem \ref{functionalrelations}. The main result of this subsection is Proposition \ref{theorem3}.
\begin{Corollary}
\label{corollaryL}
Let $s\geq 2$ be such that $s\equiv 1\pmod{q-1}.$ Then the symmetric polynomial $V_{1,s}\in K[t_1,\ldots,t_s]$ of Theorem \ref{functionalrelations}
is in fact a polynomial of $F_q[\theta][t_1,\ldots , t_s]$ of total degree $\leq s^2/(q-1)-s$ in the variables $t_1,\ldots,t_s$.
\end{Corollary}
\noindent\emph{Proof.} It follows from a simple modification of the proof of Proposition \ref{Lvalues}.
Let $\pp_1, \ldots , \pp_s$ be distinct primes in $A$, let us write 
$\aa=\pp_1\cdots\pp_s$ and let us consider the Dirichlet character $\chi$ associated to $\vartheta_{\pp_1}\cdots\vartheta_{\pp_s}$ that we also 
loosely identify with the corresponding element of $\widehat{\Delta}_\aa$. 
Since $
\aa$ is not a prime power, Remark \ref{remarktobeusedlater} implies that
\begin{equation}\label{cortobeus}
\widetilde{\pi}^{-1}L(1,\chi)g(\chi)\in \mathbb{F}_q(\chi) [\theta].\end{equation}
Now, specializing at $t_i=\zeta_i$ the root of $\pp_i$ associated to the choice of characters $\vartheta_{\pp_i}$ for all $i=1,\ldots,s$,
we obtain
$$V_{1,s}(\zeta_{1},\ldots ,\zeta_{s})=\widetilde{\pi}^{-1}L(1,\chi)g(\chi) \in\mathbb F_q(\zeta_{1}, \ldots , \zeta_{s})[\theta],$$ and the result follows from Lemma
\ref{lemma3}, the bound on the degree agreeing with that of Theorem \ref{functionalrelations}.\CVD

\subsubsection{Digit principle for the function $\omega$ and the $L$-series\label{digit}}

Let $\varphi :\TT_t\rightarrow \TT_t$ be the $\CC_\infty$-linear map defined by
$$\varphi\left(\sum_{n\geq 0}c_nt^n\right)=\sum_{n \geq 0} c_n t^{qn},\quad c_n\in\CC_\infty.$$
We also set, for $N$ a non-negative integer with its expansion in base $q$, $N=N_0+N_1q+\cdots+N_rq^r$, $N_i\in\{0,\ldots,q-1\}$:
$$\omega_N (X)=\prod _{i=0}^{r} \varphi^{i}(\omega (X))^{N_i}.$$
We then have the next Lemma.

\begin{Lemma}\label{lemmaomegavarpi}
The following identity holds:
$$\omega_N(\vartheta_\pp(\sigma_\theta)) =  \vartheta_\pp(\sigma_{\pp'})^Ng(\vartheta_\pp^N).$$
\end{Lemma}

\noindent\emph{Proof.} This is a direct application of Theorem \ref{gaussthakur}. Indeed,
$$\omega_N(\vartheta_\pp(\sigma_\theta))=\prod_{i=0}^{d-1}\omega(\vartheta_\pp(\sigma_{\theta^{q^i}}))^{N_i}=\prod_{i=0}^{d-1}
\vartheta_\pp(\sigma_{\pp'})^{q^iN_i}g(\vartheta_\pp^{q^i})^{N_i}.$$\CVD

Let $X,Y$ be two indeterminates over $K$. We introduce a family of polynomials $(G_d)_{d\geq 0}$ in $\FF_q[X,Y]$ as follows.
We set $G_0(X,Y) =1$
and
$$G_d(X,Y)= \prod_{i=0}^{d-1} (X-Y^{q^i}),\quad d\geq 1.$$ This sequence is closely related to the sequence of polynomials
$G_n(y)$ of \cite[Section 3.6]{AT}: indeed, the latter can be rewritten in terms of the former:
$$G_d(y)=G_d(T^{q^d},y^q),\quad d\geq 1,$$ in both notations of loc. cit. and ours
(\footnote{As an aside remark, we also notice that we recover in this way the coefficients of the formal series in $K[[\tau]]$ associated to Carlitz's exponential and logarithm
$$\mathfrak{e}=\sum_{i\geq 0}d_i^{-1}\tau^i,\quad\mathfrak{l}=\sum_{i\geq 0}l_i^{-1}\tau^i,$$ 
because
$d_i=G_i(\theta^{q^i},\theta)$ and $l_i=G_i(\theta,\theta^q).$ Moreover, if $\pp$ is a prime of $A$ of degree $d,$ we observe that
$$\pp=\prod_{i=1}^d(\theta-\zeta_i)=\prod_{j=0}^{d-1}(\theta-\vartheta_\pp(\sigma_{\theta^{q^j}}))=G_d(\theta, \vartheta_\pp(\sigma_\theta)).$$}). 
The polynomial $G_d$ is monic of degree $d$ in the variable $X$, and $(-1)^dG_d$ is monic in the variable $Y$ of degree 
$(q^d-1)/(q-1)$. We now define,
for $N=N_0+N_1q+\cdots+N_rq^r$ a non-negative integer expanded in base $q$, the polynomial $$\SS_{N}(t)=\prod_{i=0}^r G_i(t^{q^{i}},\theta)^{N_i}=
\prod_{i=0}^r\prod_{j=0}^{i-1}(t^{q^i}-\theta^{q^j})^{N_i}.$$
We also define the quantities associated to $N$ and $q$:
\begin{eqnarray*}
\mu_q(N)&=&\sum_{i=0}^rN_iiq^i,\\
\mu_q^*(N)&=&\frac{N}{q-1}-\frac{\ell_q(N)}{q-1},\\
\ell_q'(N)&=&\sum_{i=0}^rN_ii.
\end{eqnarray*}
\begin{Lemma}\label{propertiesSN} Let $N$ be a non-negative integer.
The following properties hold.
\begin{enumerate}
\item The polynomial $\SS_N(t)$, as a polynomial of the indeterminate $t$, is monic of degree $\mu_q(N)$.
\item As a polynomial of the indeterminate $\theta$, $\SS_N(t)$ has degree $\mu_q^*(N)$ and the leading coefficient 
is $(-1)^{\ell_q'(N)}$.
\item We have $\SS_N(\theta)=\Pi(N)$ and $v_\infty(\SS_N(\theta))=\mu_q(N)$,
where $v_\infty$ is the $\infty$-adic valuation of $\CC_\infty$.
\item We also have, for all $\zeta\in\FF_q^\text{alg}$, $v_\infty(\SS_N(\zeta))=-\mu_q^*(N)$.
\end{enumerate}
\end{Lemma}
\noindent\emph{Proof.} Easy and left to the reader.\CVD
We observe that:
$$\varphi ^d\omega (t) =\frac{1}{G_d(t^{q^{d}},\theta)}\omega(t)^{q^d}=\omega_{q^dN}(t),\quad d\geq 0$$
so that, with $N$ as above,
\begin{equation}\label{usefulident}
\omega_N(t) =\frac{\omega(t)^N }{\prod_{i=0}^r G_i(t^{q^{i}},\theta)^{N_i}}=\frac{\omega(t)^N}{\SS_N(t)}.\end{equation}

The following Proposition was inspired by a discussion with D. Goss. 
\begin{Proposition}
\label{theorem3} Let $s\geq 2$ be an integer. Let $M_1, \ldots , M_s$ be positive integers such that $M_1+\cdots + M_s \equiv 1\pmod{q-1}.$ Then:
$$W (t_1,\ldots,t_s)=\widetilde{\pi}^{-1}L(\chi_{t_1}^{M_1}\cdots \chi_{t_s}^{M_s},1)\omega_{M_1}(t_1)\cdots \omega_{M_s}(t_s)\in \mathbb{F}_q[\theta,t_1,\ldots ,t_s].$$
For all $i$, the degree in $t_i$ of $W $ satisfies
$$\deg_{t_i}(W )\leq M_i\left(\frac{\sum_jM_j}{q-1}-1\right)-\mu_q(M_i)$$
\end{Proposition}
\noindent\emph{Proof.} 
We shall write
$$H=\prod_{i=1}^s\SS_{M_i}(t_i).$$ We know from Lemma \ref{propertiesSN} that $\deg_{t_i}(H)=\mu_q(M_i)$.
Let us consider the function
$$V =\widetilde{\pi}^{-1}L(\chi_{t_1}^{M_1}\cdots \chi_{t_s}^{M_s},1)\omega^{M_1}(t_1)\cdots \omega^{M_s}(t_s),$$
so that by (\ref{usefulident}), $$V =W H.$$
Corollary \ref{corollaryL} implies that:
$$V \in \mathbb{F}_q[\theta,t_1, \cdots , t_s]$$ and we are done if we can prove that $H$ divides $V $ in $\FF_q[\theta,t_1,\ldots,t_s]$.
 
Let $\pp_1,\ldots ,\pp_s$ be distinct primes of $A$  such that $|\pp_i| -1> M_i$,
and let $\zeta_1,\ldots,\zeta_s$ be respective roots of these polynomials chosen in compatibility with 
the characters $\vartheta_{\pp_1},\ldots,\vartheta_{\pp_s}$. Let us also write
$$\chi=\vartheta_{\pp_1}^{M_1}\cdots \vartheta_{\pp_s}^{M_s}.$$ By Lemma \ref{lemmaomegavarpi},
$$\omega_{M_1}(\zeta_{1})\cdots \omega_{M_s}(\zeta_{s}) = \vartheta_{\pp_1}(\sigma_{\pp_1'})^{M_1}\cdots \vartheta_{\pp_s} (\sigma_{\pp_s'})^{M_s}g(\chi).$$
Therefore,
$$W (\zeta_1,\ldots,\zeta_s)=\widetilde{\pi}^{-1}L(1,\chi)g(\chi)\vartheta_{\pp_1}(\sigma_{\pp_1'})^{M_1}\cdots \vartheta_{\pp_s} (\sigma_{\pp_s'})^{M_s}.
$$
By (\ref{cortobeus}), 
$\widetilde{\pi}^{-1}L(1,\chi)g(\chi)\in\FF_q(\chi)[\theta]$, while $\prod_{i=1}^s\vartheta_{\pp_i}(\sigma_{\pp_i'})^{M_i}\in\FF_q(\chi)$ so that
$$W (\zeta_{1}, \cdots , \zeta_{s}) \in \FF_q(\chi)[\theta]=\mathbb F_q(\zeta_{1}, \ldots , \zeta_{s})[\theta].$$

Now, $H$ is a polynomial in $\theta$ with leading coefficient in $\mathbb F_q^\times$ (see Lemma \ref{propertiesSN}). Dividing $V $ by
$H$ as polynomials in $\theta$ we find
$$V =HQ +R,$$ where
$Q,R$ are polynomials in $\mathbb F_q[\theta, t_1, \ldots , t_s],$ and ${\rm deg}_\theta R< {\rm deg}_\theta H=\sum_i\mu_q^*(M_i)$
(the last inequality by Lemma \ref{propertiesSN}). But for $\zeta_1, \ldots, \zeta_s$ as above, we must have $Q(\theta,\zeta_1,\ldots,\zeta_s)=W (\zeta_{1}, \cdots , \zeta_{s})$ and
$$R(\zeta_{1}, \ldots , \zeta_{s})= 0.$$
This implies $R=0$ and thus $W =Q\in \mathbb F_q[\theta, t_1, \ldots , t_s].$\CVD

\subsubsection{The polynomials $W_s$}

By Proposition \ref{theorem3}, the function
$$W_s(t)=\widetilde{\pi}^{-1}L(\chi_t^s,1)\omega_N(t)=\frac{L(\chi_t^s,1)\omega(t)^s}{\widetilde{\pi}\SS_s(t)}$$
is a polynomial of $\FF_q[t,\theta]$. Furthermore, we have:
\begin{Proposition}\label{prop37}
Assuming that $s\geq2$ is an integer congruent to $1$ modulo $q-1$ and is not a power of $q$, the following properties hold.
\begin{enumerate}
\item The degree in $t$ of $W_s$ does not exceed $s(s-1)/(q-1)-s-\mu_q(s)$,
\item the degree in $\theta$ of $W_s$ is equal to $(\ell_q(s)-q)/(q-1)$.
\end{enumerate}
\end{Proposition}
By the remarks in the introduction, we know how to handle the case of $s=q^i$; we then have
$$W_{q^i}(t)=\frac{1}{\theta-t^{q^i}}.$$

\noindent\emph{Proof of Proposition \ref{prop37}.} The bound for the degree in $t$ is a simple consequence of Proposition \ref{theorem3} and Lemma \ref{propertiesSN}.
To show the property of the degree in $\theta$, we first notice that, by Lemma \ref{propertiesSN},
for all $\zeta\in\FF_q^{\text{alg}}$, 
\begin{equation}\label{valgn}
v_\infty(W_s(\zeta))=-\frac{\ell_q(s)-q}{q-1}.\end{equation}
The computation of $W_s(\zeta)$ is even explicit if $\zeta\in\FF_q$. Indeed, with the appropriate choice of a $(q-1)$-th 
root of $(\zeta-\theta)$, the fact that $\chi_\zeta=\chi_{\zeta}^s$, Lemma \ref{propertiesomega} and \cite[Theorem 1]{Pe},
\begin{eqnarray*}
W_s(\zeta)&=&
\frac{L(\chi_\zeta^s,1)\omega(\zeta)^s}{\widetilde{\pi}\SS_s(\zeta)}\\
&=&\frac{L(\chi_\zeta,1)\omega(\zeta)^s}{\widetilde{\pi}\SS_s(\zeta)}\\
&=&\frac{L(\chi_\zeta,1)\omega(\zeta)^s}{\widetilde{\pi}(\zeta-\theta)^{\frac{s-\ell_q(s)}{q-1}}}\\
&=&(\zeta-\theta)^{-\frac{1}{q-1}}(\theta-\zeta)^{-1}(\zeta-\theta)^{\frac{s}{q-1}}(\zeta-\theta)^{\frac{\ell_q(s)-s}{q-1}}
\end{eqnarray*}
and 
\begin{equation}\label{explicitvalue}
W_s(\zeta)=-(\zeta-\theta)^{\frac{\ell_q(s)-q}{q-1}}.
\end{equation}
Let us write:
$$W_s(t)=\sum_{i=0}^ga_it^i,\quad a_i\in A.$$
By (\ref{explicitvalue}), we have 
\begin{equation}\label{formulafora0}
a_0=W_s(0)=-(-\theta)^{\frac{\ell_q(s)-q}{q-1}}\end{equation}
and for all $\zeta\in\FF_q^{\text{alg}}$ we have, by (\ref{valgn}), 
$$|W_s(\zeta)|=|a_0|.$$ This means that for $i=1,\ldots,g$, $|a_i|<|a_0|$, and the identity on the degree in $\theta$ follows as well.\CVD
\begin{Corollary}\label{explicitsq}
If $\ell_q(s)=q$, then $W_s=-1$.
\end{Corollary}
\noindent\emph{Proof.} By (\ref{explicitvalue}), $W_s=a_0=-1$ in virtue of (\ref{formulafora0}).\CVD

By Corollary \ref{corollaryL}, the function $$V_{1,s}(t_1,\ldots,t_s)=\widetilde{\pi}^{-1}L(\chi_{t_1}\cdots\chi_{t_s})\omega(t_1)\cdots\omega(t_s)$$
is, for $s\equiv1\pmod{q-1}$ and $s\geq 2$, a polynomial of $A[t_1,\ldots,t_s]$. Since
$$\omega(t)=\frac{\widetilde{\pi}}{\theta-t}+o(1),$$ where $o(1)$ represents a function locally analytic at $t=\theta$,
the function $L(\chi_{t_1}\cdots\chi_{t_s},1)$ vanishes on the divisor
$$D=\bigcup_{i=1}^sD_i,$$
where $$D_i=\{(t_1,\ldots,t_{i-1},\theta,t_{i+1},\ldots,t_s)\in\CC_\infty\}.$$
In other words, in $\CC_\infty[[t_1-\theta,\ldots,t_s-\theta]]$, we have
 \begin{equation}\label{expansionLchi}
 L(\chi_{t_1}\cdots\chi_{t_s})=\sum_{i_1,\ldots,i_s\geq 1}c_{i_1,\ldots,i_s}(t_1-\theta)^{i_1}\cdots (t_s-\theta)^{i_s},\quad c_{i_1,\ldots,i_s}\in\CC_\infty,\end{equation}
 where on both sides, we have entire analytic functions (see Corollary \ref{corollaryanalytic}).
 This can also be seen, alternatively, by considering the function $F_{s-1}$ of Lemma \ref{vanishingFs}, which vanishes, and 
 observing that $$L(\chi_{t_1}\cdots\chi_{t_s},1)|_{t_i=\theta}=F_{s-1}(t_1,\ldots,t_{i-1},t_{i+1},\ldots,t_{s}).$$
 
 Let us focus on the coefficient $c_{1,\ldots,1}$ in the expansion (\ref{expansionLchi}). We then have 
 $$c_{1,\ldots,1}=\left(\left.\frac{d}{dt_1}\cdots\frac{d}{dt_s}L(\chi_{t_1}\cdots\chi_{t_s})\right)\right|_{t_1=\cdots=t_s=\theta}$$
 so that
$$V_{1,s}(\theta,\ldots,\theta)=(-1)^s\widetilde{\pi}^{s-1}\sum_{d\geq 0}\sum_{a\in A^+(d)}\frac{a'{}^s}{a}=(-1)^s\widetilde{\pi}^{s-1}c_{1,\ldots,1}\in\FF_q[\theta]$$
(by Corollary \ref{corollaryanalytic}, the series on the right-hand side is convergent). Now, by Proposition
\ref{theorem3}, $\Pi(s)$ divides the polynomial $V_{1,s}(\theta,\ldots,\theta)$ in $A$. We then set, as in the introduction:
$$\mathbb{B}_s=\frac{V_{1,s}(\theta,\ldots,\theta)}{\Pi(s)}=G_s(\theta)\in A.$$
\subsection{Proof of Theorem \ref{bernoulli}}

We begin the proof with a couple of simple remarks.
Firstly, if $B$ is a polynomial of $A[t]$ and if $\pp$ is a prime of degree $d>0$, then 
$$
\tau^dB\equiv B\pmod{\pp}.
$$
The reason for this is that $\pp$ divides the polynomial $\theta^{q^d}-\theta$. In particular,
\begin{equation}\label{congruence}
(\tau^dB)(\theta)\equiv B(\theta)\pmod{\pp}.
\end{equation}
Secondly, recalling the $\CC_\infty$-linear operator $\varphi$ of subsection \ref{digit}, 
we have 
$$\tau\varphi=\varphi\tau=\rho,$$ where $\rho$ is the operator defined by $\rho(x)=x^q$ for all $x\in\CC_\infty((t))$.
In particular, if $s=\sum_{i=0}^rs_iq^i$ is expanded in base $q$ and if $d\geq r\geq i$, from
$$\tau^d\varphi^i=\tau^{d-i}\tau^i\varphi^i=\tau^{d-i}\rho^i$$
we deduce
$$
(\tau^d\omega_s)(t)=\prod_{i=0}^r((\tau^{d-i}\omega)(t))^{s_iq^i},
$$
so that 
\begin{equation}\label{secrem}
(\tau^d\omega_s)(t)=\prod_{i=0}^rG_{d-i}(t,\theta)^{s_iq^i}\omega(t)^s.
\end{equation}
We can finish the proof of Theorem \ref{bernoulli}.
By (\ref{congruence}), 
$$\mathbb{B}_s\equiv (\tau^dW_s)(\theta).$$ We shall now compute $(\tau^dW_s)(\theta)$.
If $d>r$, we can write
$$
G_{d-i}(t,\theta)^{s_iq^i}=(t-\theta)^{s_iq^i}\prod_{j=1}^{d-i-1}(t-\theta^{q^j})^{s_iq^i},
$$
and 
$$
\prod_{i=0}^rG_{d-i}(t,\theta)^{s_iq^i}=(t-\theta)^sF(t),
$$
where $F(t)$ is a polynomial such that
$$
F(\theta)=\prod_{i=0}^rl_{d-i-1}^{s_iq^i}.
$$
Since
$$
(\tau^dW_s)(t)=\widetilde{\pi}^{-q^d}L(\chi_t^s,q^d)(t-\theta)^s\omega(t)^sF(t)
$$
and $\lim_{t\rightarrow\theta}(t-\theta)\omega(t)=-\widetilde{\pi}$,
we get
\begin{eqnarray*}
\lim_{t\rightarrow\theta}(\tau^dW_s)(t)&=&(-1)^s\widetilde{\pi}^{-q^d}\zeta(q^d-s)\widetilde{\pi}^s\prod_{i=0}^rl_{d-i-1}^{s_iq^i}\\
&=&(-1)^s\frac{BC_{q^d-s}}{\Pi(q^d-s)}\prod_{i=0}^rl_{d-i-1}^{s_iq^i}.
\end{eqnarray*}
Our Theorem \ref{bernoulli} follows at once.\CVD
\begin{small}
\noindent\emph{Acknlowledgements.} The authors would like to warmly thank Vincent Bosser, David Goss, Matthew Papanikolas, Rudolph Perkins, Denis Simon, Lenny Taelman, Floric Tavares-Ribeiro, Dinesh Thakur for fruitful discussions about the themes developed in this 
work. In addition, we are indebted with: David Goss who suggested the appropriate direction of investigation in order to obtain Proposition 
\ref{theorem3}, Denis Simon who discovered Lemma \ref{simon}, and Lenny Taelman who performed 
numerical computations on Bernoulli-Carlitz fractions providing additional evidence to our Conjecture \ref{conj1}. Moreover, 
we would like to thank the {\em Istituto de Giorgi} of Pisa, Italy, for having hosted the workshop 
``Zeta functions and $L$-series in positive characteristic" in November 2012. We have participated to this workshop and 
this allowed us to fruitfully discuss with colleagues at
the same time we were involved in the final editing of this work, in the very nice conditions the {\em Istituto} offered.
\end{small}

\end{document}